\documentclass{siamart1116}

\usepackage{hyperref}
\usepackage{bm}
\usepackage{float}
\usepackage{amsmath}
\usepackage{amsfonts}
\usepackage{graphicx}
\usepackage{xcolor}
\usepackage{comment}
\usepackage{framed}

\newcommand{\mcF}{\mathcal{F}}

\newcommand{\mcI}{\mathcal{I}}

\newcommand{\mcX}{\mathcal{X}}

\newcommand{\mcL}{\mathcal{L}}

\newcommand{\mcA}{\mathcal{A}}
\newcommand{\mcT}{\mathcal{T}}

\newcommand{\gam}{{\gamma}}

\def \xb{\bm{x}}
\def \yb{\bm{y}}

\def \sl{{s_l}}

\def\RRR{{\Re}}
\def\TTT{{T}}
\def\vecbf{\vec}
\def\vecupn{{({\vecbf u}_n)}}
\def\vecgpn{{({\vecbf g}_n)}}

\def\vecbpn{{({\vecbf b}_n)}}
\def\vecun{{\vecbf u}_n}

\def\vecgn{{\vecbf g}_n}

\def\vecbn{{\vecbf b}_n}

\def\vecup{{({\vecbf u})}}
\def\vecgp{{({\vecbf g})}}

\def\vecu{{\vecbf u}}

\def\vecg{{\vecbf g}}

\def\vecb{{\vecbf b}}

\def\hangone{\vskip0pt\hangindent=31pt\hangafter=1}
\def\hangzero{\vskip0pt\hangindent=19pt\hangafter=1}
\def\hangtwo{\vskip0pt\hangindent=60pt\hangafter=1}

\usepackage{color}

\definecolor{officegreen}{rgb}{0.0, 0.5, 0.0}

\newsiamthm{remark}{Remark}

\title{A general framework for substructuring-based domain decomposition methods for\\ models having nonlocal interactions}
\author{
Giacomo Capodaglio\thanks{Department of Scientific Computing, Florida State University, Tallahassee FL 32306; Current address: Computational Physics and Methods, Los Alamos National Laboratory, Los Alamos NM 87545;
{\tt gcapodaglio@lanl.gov}.
}
\and
Marta D'Elia\thanks{Computational Science and Analysis, Sandia National Laboratories, Livermore CA 94550; 
{\tt mdelia@sandia.gov}
}
\and
Max Gunzburger\thanks{Department of Scientific Computing, Florida State University, Tallahassee FL 32306;
{\tt mgunzburger@fsu.edu}.
}
\and
Pavel Bochev\thanks{Center for Computing Research, Sandia National Laboratories, Albuquerque NM 87321;
{\tt pbboche@sandia.gov}
}
\and
Manuel Klar\thanks{Department of Mathematics, Universit{\"a}t Trier, 54296 Trier, Germany; {\tt klar@uni-trier.de}
}
\and
Christian Vollmann\thanks{Department of Mathematics, Universit{\"a}t Trier, 54296 Trier, Germany; {\tt vollmann@uni-trier.de}
}
}

\begin{document}
\maketitle

\begin{abstract}
A rigorous mathematical framework is provided for a substructuring-based domain-decomposition approach for nonlocal problems that feature  interactions between points separated by a finite distance. Here, by  \emph{substructuring} it is meant that a traditional geometric configuration for local partial differential equation problems is used in which a computational domain is subdivided into non-overlapping subdomains. In the nonlocal setting, this approach is \emph{substructuring-based} in the sense that those subdomains interact with neighboring domains over interface regions having finite volume, in contrast to the local PDE setting in which  interfaces are lower dimensional manifolds separating abutting subdomains. Key results include the equivalence between the global, single-domain nonlocal problem and its multi-domain reformulation, both at the continuous and discrete levels. These results provide the rigorous foundation necessary for the development of efficient solution strategies for nonlocal domain-decomposition methods.
\end{abstract} 


%nonlocal, FEM, integral operators, integral equations, numerical analysis

\pagestyle{myheadings}
\thispagestyle{plain}
\markboth{G. Capodaglio, M. D'Elia, M. Gunzburger, P. Bochev, M. Klar and C. Vollman}
{{A General Framework for Nonlocal Domain Decomposition}}

%%%%%%%%%%%%%%%%%%%%%%%%%%%%%%%%%%%%%%%%%%%%%%%%%%%%%%%%%%%%%%%%%%%%%%%%
%%%%%%%%%%%%%%%%%%%%%%%%%%%%%%%%%%%%%%%%%%%%%%%%%%%%%%%%%%%%%%%%%%%%%%%%
\section{Introduction}\label{intro}

Nonlocal models have become a popular alternative to partial differential equation (PDE) models due to their ability to describe effects that PDEs fail to capture. In particular, a nonlocal model can describe multiscale and anomalous behavior for applications that exhibit hierarchical features that cannot be reproduced by a classical model. These applications include, among others, subsurface transport \cite{Benson2000,Schumer2003,Schumer2001},
image processing \cite{Buades2010,DElia2019imaging,Gilboa2007,Lou2010},
multiscale and multiphysics systems \cite{Alali2012,Askari2008},
magnetohydrodynamic \cite{Schekochihin2008},
finance \cite{Scalas2000,Sabatelli2002}, and stochastic processes \cite{Burch2014,DElia2017,Meerschaert2012,Metzler2000,Metzler2004}.

The general class of nonlocal models we consider are characterized by integral operators having the form
\begin{equation}\label{eq:simple-laplacian}
\mcL u(\xb)= \int_{B_\delta(\xb)} (u(\yb)-u(\xb))\gam(\xb,\yb)\,d\yb,
\end{equation}
where $B_\delta(\xb)$ denotes the ball (usually Euclidean) centered at $\xb$ with radius $\delta$ (usually referred to as the {\em horizon} or {\em interaction radius}) and $\gam(\xb,\yb)$ is an application-dependent kernel function (usually symmetric in its arguments and nonnegative) that determines the regularity properties of the solution. The nonlocality inherent in \eqref{eq:simple-laplacian} is clear: points $\xb$ interact with all points $\yb$ located within a distance $\delta$. Compared to that for the local PDE setting, the integral form clearly reduces regularity requirements on the solution and allows for the capture of long-range interactions. 

However, the utilization of nonlocal models in applications that {could benefit from their} improved predictive capabilities is hindered by several modeling and numerical challenges. {These include} the unresolved treatment of nonlocal interfaces \cite{Alali2015,Capodaglio2020}, the nontrivial prescription of nonlocal volume constraints (the nonlocal counterpart of boundary conditions) \cite{Cortazar2008,DEliaNeumann2019}, and the fact that computational costs attendant to the use of nonlocal problems may become prohibitive as the extent of the nonlocal interactions increases; see, e.g., \cite{acta20,DEliaFEM2020} for variational methods and \cite[Chapter 7]{Chen2006meshless} for mesh-free methods. Other critical challenges are related to the uncertain nature of model parameters; in fact, modeling parameters such as $\delta$ and those characterizing the kernel, applied forces, and/or sources can be non-measurable, sparse, and/or subject to noise. Research on such topics is very active (see, e.g., \cite{Antil2019FracControl,Antil2018sSpectralControl,DElia2019FracControl,DElia2014DistControl,DElia2016ParamControl,DElia2019imaging,Gulian2019,Pang2020,Pang2019fPINNs,Pang2017discovery,You2020Regression}) but further consideration of them is beyond the scope of this work. 

Here, we focus on the treatment of nonlocal interfaces and, notably on the design of nonlocal domain-decomposition (DD) formulations with the aim of reducing computational costs by increasing the parallel concurrency in the numerical solution of nonlocal problems. Specifically, the goal is to address the high computational cost associated with nonlocal models by providing a foundational algorithmic framework for their parallel solution, mirroring that of successful parallel {DD} algorithms for PDEs such as, e.g., Finite Element Tearing and Interconnecting (FETI) \cite{farhat1991method} and other approaches \cite{mathew2008domain,toselli2006domain}.

This work is part of a comprehensive effort by the authors to fill the theoretical and practical gaps in the current understanding of nonlocal interfaces {(both physical ones and those created by DD solution algorithms)} by developing a rigorous nonlocal interface theory for nonlocal diffusion (see the preliminary work \cite{Capodaglio2020}), including pure fractional diffusion, and nonlocal mechanics. Our ultimate goal is to design efficient and scalable DD solvers to unlock the full potential of nonlocal models. To this end, as is often done in nonlocal modeling, we {draw inspiration from the existing work on} classical DD methods for PDEs. 
{Notwithstanding the vast knowledge accumulated by the DD community}, extension of local DD methods to nonlocal models {remains} a nontrivial task: nonlocality introduces many challenges and limitations.
Some of these challenges are shared with local DD methods. An example of a shared challenge is the proper treatment of {\it floating subdomains}\footnote{In a multi-domain setting, floating domains are subdomains of the decomposition that are either internal (they do not share boundaries with the physical boundary) or share boundaries only with parts of the physical boundary at which Neumann-type conditions are prescribed. In this context ``floating'' refers to the fact that these domains either do not have volume constraints because they are internal, or have partial constraints of Neumann type.} which require special attention due to the singularity of the discretized equations on such subdomains.

However, other challenges are unique to the nonlocal setting and require new approaches that have no {analogues} in the local setting. 
For example, in the local PDE setting, {many} DD methods {typically start} by breaking {up} the computational domain into non-overlapping subdomains, a process we refer to as \emph{substructuring of the domain}. These subdomains interact only through their shared boundaries on which one usually imposes some appropriate continuity conditions. 
Although one starts from the same initial geometric configuration in which the domain is substructured into non-overlapping subdomains, inherent nonlocal interactions between the subdomains force one to expand these subdomains to include parts of neighboring subdomains having nonzero volume, causing an overlap of a thickness determined by the interaction radius $\delta$. This overlap is required because, in the nonlocal setting, it is not possible to define subdomain problems simply by restricting the global operator to the subdomains. For this reason we refer to our approach as being \emph{substructuring-based}. 
It is important to point out that the thickness of the overlap regions depends solely on the {\em modeling} parameter $\delta$ and is unrelated to the discretization method employed and, in particular, to the grid size. 
As a result, discretization of the decomposed nonlocal problem requires special care because the overlapping regions induced by the decomposition do not, in general, match the underlying mesh.

The current state of the art of nonlocal DD methods is very limited, with \cite{aksoylu2011variational} being perhaps the most relevant work. In that paper, the authors consider a simple two-domain configuration and develop a variational approach to DD based on adding an interface equation and a new variable that lives on the overlap between the subdomains, while using test functions that vanish on the interface for each subdomain. Decoupling is achieved by solving a Schur-complement problem for the interface variable, similar to a conventional FETI scheme. The subproblem definition in \cite{aksoylu2011variational} does not consider  multi-domain configurations nor does it consider floating subdomains. 
{As a result, extension of this approach to the general multi-domain case, which is a prerequisite for efficient parallel nonlocal DD algorithms, remains an open question.}

In this paper we formulate a general framework for nonlocal DD problems {that has the following equivalence property:}
\begin{equation}\label{goalgoal2}
\begin{aligned}
&\mbox{\em the discrete solution obtained via the DD approach is identical}
\\[-.75ex]
&\mbox{\em  to the discrete solution obtained for the parent single domain.}
\end{aligned}
\end{equation} 
{While the significance of this property is self-evident, its fulfillment is far from trivial, and is one of the key challenges addressed in the paper.} 

As already mentioned, we refer to our approach as ``substructuring-based'' because, much like as it is in standard non-overlapping DD, the subdomains interact only through their shared interfaces. Of course, the key difference is that in the nonlocal setting these shared interfaces are nonlocal, i.e., they  are  regions having finite volumes as opposed to the local case in which interfaces are lower-dimensional manifolds.

The main contributions of this paper are as follows.
\begin{itemize}
    \item We introduce a systematic way to decompose the domain given an existing mesh for the single-domain problem and we discuss ways to make the decomposed domains compatible with the given mesh. Specifically, we provide a recipe for decomposing the domain that prevents integration over partial (cut) finite elements by using approximate neighborhoods in a manner such {that the equivalence of the decomposed solution and the single-problem one is not compromised.}
    \item We formulate a continuous nonlocal DD system of subdomain problems and prove that it is equivalent to the single-domain problem, {i.e., we provide a solution for \eqref{goalgoal2}.} The key ingredient is the appropriate definition of indicator functions that keep track of the number of overlapping subdomains. 
    \item {We define a Galerkin finite element discretization of the nonlocal DD problem and prove that it is equivalent to the discretization of the single-domain domain problem effected using the same type of finite element functions. This equivalence holds for the finite-dimensional variational formulation and for the corresponding matrix form.}
\end{itemize}

Our nonlocal DD formulation provides a mathematical foundation for the development of a range of efficient numerical algorithms for the parallel solution of nonlocal problems that mirror existing approaches for local problems. For example, treating the nonlocal interface equations as constraints and using Lagrange multipliers to enforce them lends itself to the development of nonlocal FETI \cite{farhat1991method} or Arlequin-like \cite{BenDhia_05_IJNME}  algorithms.  Alternatively, one can choose to view these coupling conditions as an optimization objective and treat the subdomain equations as constraints. Such an approach would lead to nonlocal optimization-based DD methods that are nonlocal counterparts of the methods in \cite{Gunzburger_00_AMC,Gunzburger_99_CMA,Bochev_17a_SISC}.

It should be noted though that realizing the potential of our DD framework to reduce the computational burden of solving nonlocal problems requires the ratio between $\delta$ and the diameter of the subdomains to be smaller than 1, i.e., we target problems for which the extent of the nonlocal interactions is much smaller than the diameter of the domain. Such problems arise in several engineering applications such as, e.g., nonlocal mechanics, and are the main motivator for this work. In contrast, for applications described by nonlocal operators with infinite interactions, a DD approach may not be as effective because  the interaction regions would span a portion of the domain that is of the same size (or even larger) than the domain itself.

Finally, we mention that often one may be given a decomposition of $\Omega$ into a {\em few} subdomains which {are} constructed to follow well-defined geometric entities, e.g., a wing and a fuselage, or different media properties, e.g., different diffusion coefficients, within $\Omega$. Such ``physically''-motivated DDs typically arise in the context of \emph{mesh tying} \cite{Laursen_03_IJNME,Bochev_07_CMAME} in which a complex geometric entity is broken into smaller parts to enable efficient mesh generation. In contrast, here we focus on DD as a means for faster and more efficient parallel solution methods for nonlocal problems in which case the number of the subdomains is very large and they do not generally follow any ``physics''-motivated interfaces.

%%%%%%%%%%%
\smallskip\paragraph{Outline of the paper} The paper is organized as follows. In Section \ref{sec:vcp}, we recall the variational formulation of a single-domain volume-constrained nonlocal Poisson problem and briefly describe its discretization via finite element methods. In Section \ref{sec:continuousDD}, we introduce the continuous formulation of a multi-domain {nonlocal} DD method and prove its equivalence to the single-domain problem presented in Section \ref{sec:vcp}. Section \ref{femsubsys} explains how we address the decomposition problem by formulating rules for the construction of the subdomains and their interaction regions that {fulfill \eqref{goalgoal2}, i.e.,}  {guarantee the equivalence of the DD and the single-domain problems.} In the same section we also introduce the discretized subproblems and their matrix forms and {show their equivalence  to the underlying single-domain formulation.} Concluding remarks are provided in Section \ref{sec:conclusion}. 

{In the technical report  \cite{sandiareport}, we further elucidate the equivalence between the multi-domain formulation and the single-domain problem and illustrate the application of the nonlocal DD framework.  Specifically, we use the framework developed in this work along with a FETI solution approach to obtain, for a very simplified setting, illustrative numerical examples of nonlocal DD problems.}

%%%%%%%%%%%
\section{A nonlocal (single-domain) volume-constrained problem and its finite element discretization}\label{sec:vcp}

\begin{figure}[h!]
   \centering
   \includegraphics[width=2.2in]{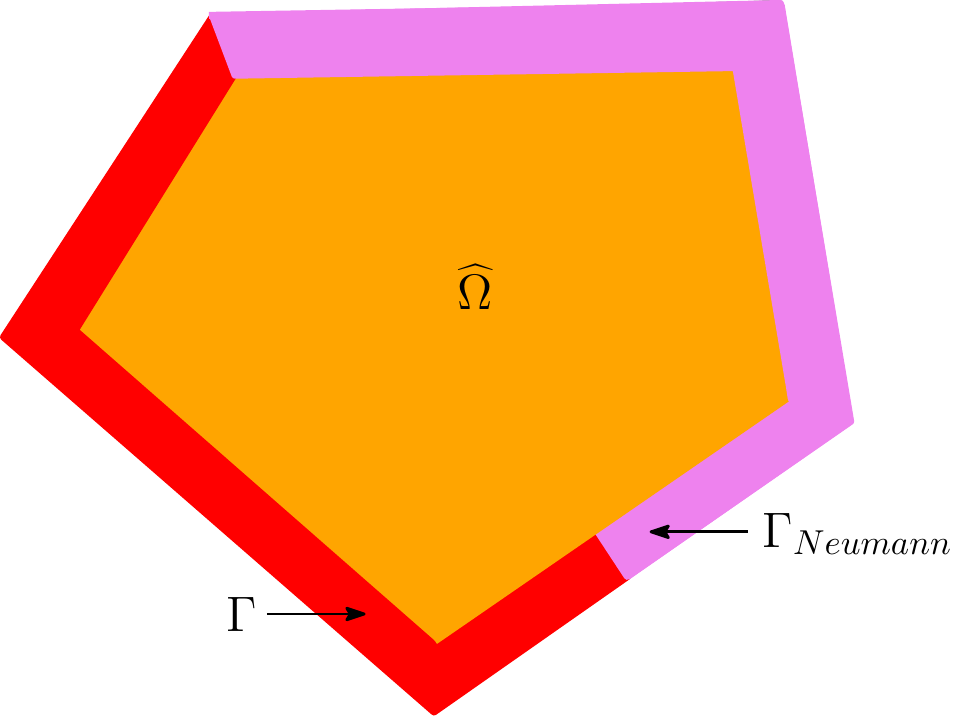}
   \quad
   \includegraphics[width=2.in]{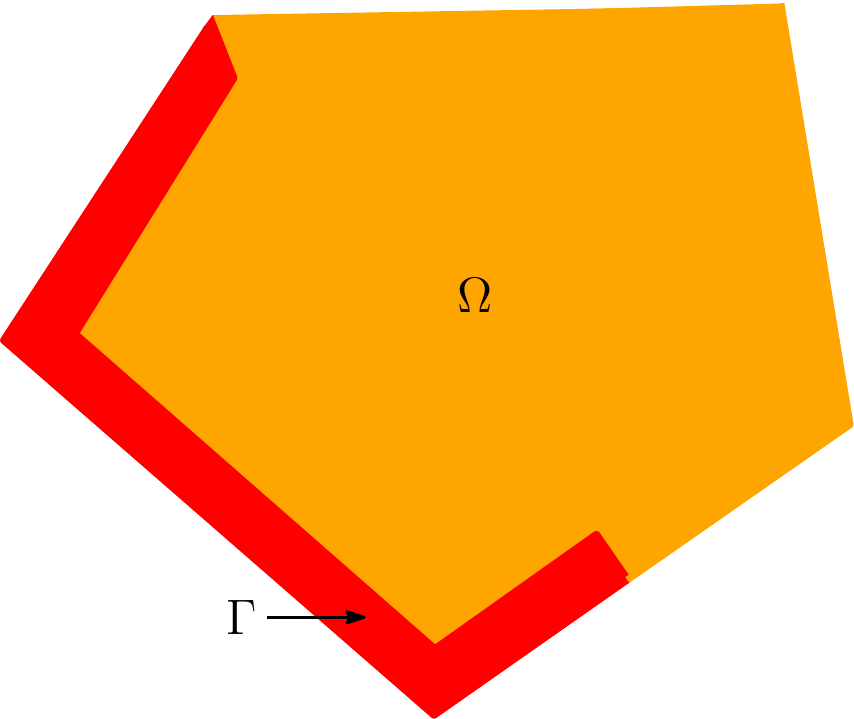} 
\caption{Left: A domain $\widehat\Omega$ and its associated interaction domain $\Gamma\cup\Gamma_{N\!eumann}$ on which Dirichlet and Neumann volume constraints are imposed on $\Gamma$ and $\Gamma_{N\!eumann}$, respectively. Right: the domain $\Omega=\widehat\Omega\cup\Gamma_{N\!eumann}$.}
   \label{fig1.1}
\end{figure}

For simplicity, in this work, we consider the two-dimensional case. Let $\widehat\Omega$ denote a bounded, open subset of $\RRR^2$. For any $\delta > 0$, often referred to as the {\em horizon} or {\em interaction radius}, we define the associated {\em interaction domain} as the closed region
{\begin{equation}\label{eq:int_dom}
\Gamma_{interaction}= \{\yb \in  \RRR^2  \setminus \widehat\Omega\,:\, \exists \;\xb \in \widehat\Omega
\;\; \hbox{such that} \;\; |\xb - \yb| \leq \delta \; \} \,.
\end{equation}}
The interaction domain is split into two disjoint parts $\Gamma$ and $\Gamma_{N\!eumann}$, where $\Gamma$ is a nonempty closed domain, whereas $\Gamma_{N\!eumann}$ is allowed to be empty. Thus, we have that $\Gamma\cup\Gamma_{N\!eumann}=\Gamma_{interaction}$ and $\Gamma\cap\Gamma_{N\!eumann}=\emptyset$, where $\Gamma_{N\!eumann}$ is open along its common boundary with $\Gamma$. Also, note that $\Gamma_{interaction}$ and therefore also $\Gamma$ and $\Gamma_{N\!eumann}$ depend on $\delta$, even though that dependence is not explicitly indicated. Figure \ref{fig1.1}-left illustrates this geometric configuration.\footnote{Domains such as $\Gamma$ and $\Gamma_{N\!eumann}$ in Figure \ref{fig1.1}-left and therefore also in subsequent figures are stylized versions of their true shapes. For example, because points $\xb$ interact only with points $\yb\in B_\delta(\xb)$, where $B_\delta(\xb)$ denotes the Euclidean ball of radius $\delta$ centered at $\xb$, those domains have rounded corners. However, in practice, one can keep the stylized domains because the points outside the true interaction domains are not accessed during a finite element assembly process.} 
%We assume throughout that $\Gamma\ne\emptyset$ whereas $\Gamma_{N\!eumann}$ is allowed to be empty.
%%%

The strong formulation of a nonlocal volume-constrained Poisson problem is given by \cite{acta20,du19,du12,du13}\footnote{The problem \eqref{eq:strong_nonloc1} is a nonlocal analogue of the PDE Poisson problem $-\nabla\cdot(\kappa\nabla u)=f_{
\widehat\Omega}$ in $\widehat\Omega$, $u=g$ on a nonempty part of the boundary of $\widehat\Omega$, and $\kappa\nabla u\cdot {\bf n}= f_{N\!eumann}$ on the rest of that boundary.}
\begin{equation}\label{eq:strong_nonloc1}
\begin{cases}
\displaystyle  
 -2\int_{\widehat\Omega\cup\Gamma\cup\Gamma_{N\!eumann}} 
\big(u(\yb) -u(\xb)\big)\gamma(\xb,\yb)d\yb 
= f_{\widehat\Omega}(\xb) & \xb\in\widehat\Omega
\qquad\qquad\, (a)
\\[4mm]
u(\xb) = g(\xb) &  \xb\in\Gamma \qquad\qquad\,\, (b)
\\[2mm]
 2\displaystyle\int_{\widehat\Omega\cup\Gamma\cup\Gamma_{N\!eumann}} 
\big(u(\yb) -u(\xb)\big)\gamma(\xb,\yb)d\yb 
= f_{N\!eumann}(\xb) & \xb\in\Gamma_{N\!eumann}, \,\, (c)
\end{cases}
\end{equation}
where $f_{\widehat\Omega}(\xb)$, $f_{N\!eumann}(\xb)$, and $g(\xb)$ are given functions and $\gamma(\xb,\yb)$ is a  given symmetric nonnegative kernel, i.e., $\gamma(\xb,\yb) = \gamma(\yb,\xb)$ for $\xb,\,\yb\in \widehat\Omega\cup\Gamma\cup\Gamma_{N\!eumann}$.
{Equations (\ref{eq:strong_nonloc1}b) and (\ref{eq:strong_nonloc1}c) are nonlocal analogues of Dirichlet and Neumann boundary conditions, respectively, for PDEs; specifically, equation (\ref{eq:strong_nonloc1}b) is a {\it Dirichlet volume constraint} imposed on a set $\Gamma$ having nonzero measure in $\Re^2$ and, if $\Gamma_{N\!eumann}$ also has nonzero measure in $\Re^2$, (\ref{eq:strong_nonloc1}c) is a {\it Neumann volume constraint} on that set.}

The nonlocal operators in (\ref{eq:strong_nonloc1}a) and (\ref{eq:strong_nonloc1}c) are identical up to a sign so that these equations can be combined to obtain an equivalent, more compact,  strong form
\begin{equation}\label{eq:strong_nonloc}
\begin{cases}
\displaystyle  
-2\int_{\Omega\cup\Gamma} 
\big(u(\yb) -u(\xb)\big)\gamma(\xb,\yb)d\yb 
= f(\xb) & \xb\in\Omega \\[2.5mm]
u(\xb) = g(\xb) &  \xb\in\Gamma,
\end{cases}
\end{equation}
where $\Omega=\widehat\Omega\cup\Gamma_{N\!eumann}$, $f(\xb)|_{\widehat\Omega}=f_{\widehat\Omega}(\xb)$, and $f(\xb)|_{\Gamma_{N\!eumann}}=-f_{N\!eumann}(\xb)$. The strong form \eqref{eq:strong_nonloc} corresponds to the configuration in Figure \ref{fig1.1}-right and is used in the remainder of the paper.
{We assume that $\Omega$ does not include the boundary portion $\partial\Omega\cap\partial\Gamma$, but does include the boundary portion $\partial\Omega\setminus(\partial\Omega\cap\partial\Gamma)$, where $\partial\Omega$ and $\partial\Gamma$ denote the boundaries of $\Omega$ and $\Gamma$, respectively.}

We define the function spaces
\begin{equation}\label{globalspace}
\left\{\begin{aligned}
& W= \{w \in L^2(\Omega\cup\Gamma) \,\,:\,\, |||w||| < \infty\} \\
   &\qquad\mbox{where} \quad 
  |||w|||^2 =  \int_{\Omega\cup\Gamma}\int_{\Omega\cup\Gamma}
  |w(\yb) - w(\xb)|^2\gamma(\xb,\yb)d\yb  d\xb +
  \|w\|^2_{L^2(\Omega\cup\Gamma)}\\
& W^0 = \{w\in W\,\,:\,\, w=0 \,\,\mbox{for $\xb\in\Gamma$}\}
\end{aligned}\right.
\end{equation}
and, for $u,v\in W$, we define the bilinear form $\mathcal{A}(\cdot,\cdot)$ and linear functional $\mcF(\cdot)$ as
\begin{equation}\label{eq:bil1}
\left\{ 
\begin{aligned}
&\mathcal{A}(u,v) = 
\int_{\Omega\cup\Gamma}  \int_{\Omega\cup\Gamma}\big(u(\yb) - u(\xb)\big)\big(v(\yb) - v(\xb)\big)\gamma(\xb,\yb)d\yb  d\xb 
\\
&\mcF(v) = \int_\Omega v(\xb)f(\xb)d\xb. 
\end{aligned}\right.
\end{equation}
Then, a weak formulation of the nonlocal volume-constrained problem  \eqref{eq:strong_nonloc} can be stated as \cite{acta20,du19,du12,du13}
\begin{equation}\label{eq:weak_nonloc}
\begin{aligned}
&\mbox{\em given $f(\xb)\in W'$, $g(\xb)\in W_\Gamma$, and a kernel $\gamma(\xb,\yb)$, find $u(\xb) \in W$ such that}
\\
&  \qquad\qquad 
    \mathcal{A}(u,v) = \mcF(v) \quad \forall\, v \in W^0
    \qquad\mbox{\em subject to $u(\xb)=g(\xb)$ for $\xb\in\Gamma$}. 
\end{aligned}
\end{equation}
Here, $W'$ denotes the dual space of bounded linear functionals on $W^0$ with respect to the standard $L^2$ duality pairing and $W_\Gamma$ denotes the nonlocal ``trace'' space defined as $W_\Gamma=\{w|_\Gamma\,\,:\,\, w\in W\}$.
{If $\Gamma$ has nonzero measure in $\Re^2$,} the well posedness of the problem \eqref{eq:weak_nonloc} is proved in, e.g., \cite{acta20,du19,du12,du13}.   

For $v\in W$, define the energy functional
\begin{equation}\label{singleenergy}
\mathcal{E}_{single}(v) = \frac{1}{2}\mathcal{A}(v,v) - \mcF(v).
\end{equation}
Then, \eqref{eq:weak_nonloc} is equivalent to the  minimization problem \cite{acta20,du19,du12,du13}
\begin{equation}\label{eq:minimiz_nl}
\begin{aligned}
&\mbox{\em given $f(\xb)\in W'$, $g(\xb)\in W_\Gamma$, and a kernel $\gamma(\xb,\yb)$, find $u(\xb) \in W$ such that}
\\
&  \qquad\qquad
     \mathcal{E}_{single}(u) = \inf\limits_{v\in W} {\mathcal E}_{single}(v) 
     \qquad\mbox{\em subject to $u(\xb)=g(\xb)$ for $\xb\in\Gamma$}.
\end{aligned}
\end{equation}

\begin{remark}\label{kernelchoice}
{\em 
The functional setting and the well posedness of the nonlocal problem depend on the kernel  $\gamma(\xb,\yb)$. For example, if the kernel is square integrable (i.e., $\int_{\Omega\cup\Gamma}\big(\gamma(\xb,\yb)\big)^2 d\yb<\infty$ for all $\xb\in\Omega\cup\Gamma$) or if the kernel is integrable and translationally invariant (i.e., $\int_{\Omega\cup\Gamma}\gamma(\xb,\yb) d\yb<\infty$ for all $\xb\in\Omega\cup\Gamma$ and $\gamma(\xb,\yb))=\gamma(\yb-\xb)$), it is known that $W=L^2(\Omega\cup\Gamma)$; see, e.g., \cite{acta20,du19,du12,du13}. On the other hand, for the fractional kernel $\gamma(\xb,\yb)\propto |\yb-\xb|^{-d-2s}$ with $d$ denoting the space dimension and $0<s<1$, it is known that $W=H^s(\Omega\cup\Gamma)$, i.e., a fractional Sobolev space; again, see, e.g., \cite{acta20,du19,du12,du13}.
However, the technical aspects of this work are largely independent of the choice of the kernel as long as the nonlocal problem remains well posed. Moreover, some conditions such as the symmetry of  $\gamma(\xb,\yb)$ can be further relaxed under some additional assumptions \cite{acta20,DElia2017} that ensure the well posedness of the nonlocal problem. Likewise, one can also relax the condition that $\gamma(\xb,\yb)$ is nonnegative everywhere; see, e.g., \cite{medu}. {Although for such kernels, an energy minimization characterization of the nonlocal problem may not be available, 
well-defined weak formulations of the strong form nonlocal equations still exist {for more general cases}. As a result,} the algorithms developed in this work can be extended in a straightforward manner to problems that cannot be characterized in terms of an energy minimization setting.
}\qquad$\Box$
\end{remark}

%%%%%%%%%%%%%%%%%%%
\subsection{Finite element discretization of the nonlocal volume-constrained problem}\label{sec1-1}
The weak formulation \eqref{eq:weak_nonloc} of the nonlocal volume-constrained problem can be discretized using a finite element method as follows. Let $\mathcal{T}^h$ denote a finite element triangulation of $\Omega \cup \Gamma$ parameterized by a grid size parameter $h$. We assume that $\mathcal{T}^h$ conforms to the boundary of $\Omega$, i.e., $\partial\Omega$ consists of finite element edges. This requirement can be satisfied by first constructing a grid in $\Omega$ after which a grid is constructed in $\Gamma$ that shares element vertices with those of the grid in $\Omega$ along their common boundary. For simplicity, in the sequel, we restrict the discussion to Lagrangian finite element spaces.

Let $W^h\subset W$ and {$W^{0,h}\subset W^0$} denote finite element (FE) subspaces. Then, a FE approximation $u^h(\xb)\in W^h$ of the solution $u\in W$ of \eqref{eq:weak_nonloc} is defined to be the solution of the discretized weak formulation
\begin{equation}\label{eq:weak_nonlochh}
\begin{aligned}
&\mbox{\em given $f(\xb)\in W'$, $g(\xb)\in W_\Gamma$, and a kernel $\gamma(\xb,\yb)$}, 
\\& \mbox{\em find $u^h(\xb) \in W^h$ such that}
\\
&  \qquad 
    \mathcal{A}(u^h,v^h) = \mcF(v^h) \quad \forall\, v^h \in W^{0,h}\qquad
\mbox{\em subject to $u^h(\xb)=g^h(\xb)$ for $\xb\in\Gamma$}, 
\end{aligned}
\end{equation}
where $g^h(\xb)$ denotes an approximation of $g(\xb)$ that is usually chosen to be the FE interpolant\footnote{When $g(\xb)$ is not of class $C^0$, $g^h(\xb)$ can be defined by, e.g., least-squares approximation or Clement interpolation \cite{Ciarlet_02_BOOK}.} of  $g(\xb)$. As long as {$\Gamma$ has nonzero measure in $\Re^2$,} the well posedness of problem \eqref{eq:weak_nonlochh} is guaranteed; see, e.g., \cite{acta20,du19,du12,du13}. 

Let $\widetilde N^h$ denote the number of degrees of freedom corresponding to the nodes in $\Omega\cup\Gamma$ and let $N_{h}$ denote the number of degrees of freedom corresponding to the (possibly semi-) open domain $\Omega$, with the remaining and $\widetilde N^h-N^h$ degrees of freedom corresponding to the closed domain $\Gamma$.
Note that because $\Omega$ is an open domain with respect to its common boundary with $\Gamma$ and $\Gamma$ itself is a closed domain, nodes and degrees of freedom along their common boundary are assigned to $\Gamma$. We then define the finite element subspace $W^h\subset W$ as the span of a nodal finite element basis $\{\phi_i(\xb)\}_{i=1}^{\widetilde N^h}$ so that a finite element approximation $u^h(\xb)$ of the solution $u(\xb)$ of \eqref{eq:weak_nonloc} can be expressed as
\begin{align}\label{eq:fem_soln}
u^h(\xb) = \sum_{i=1}^{N^h} \vecup_i \phi_i(\xb)
+ \sum_{i=N^h+1}^{\widetilde N^h}\vecgp_i \phi_i(\xb),
\end{align}
where $\vecu$ denotes an $N^h$-vector of unknown coefficients and $\vecg$ denotes an $(\widetilde N^h-N^h)$-vector of nodal values of the approximation $g^h(\xb)$ of $g(\xb)$. 
We recall that the support of each nodal basis function $\phi_i(\xb)$ comprises all elements sharing the node $\xb_i$. As a result, all basis functions corresponding to nodes in $\Omega$ vanish on $\Gamma$ and  span$\{\phi_i(\xb)\}_{i=1}^{N^h}\subset W^{0,h}$. 

Let $\mathbb{A}_{single}$  and $\vecbf{b}_{single}$ denote the $N^h\times N^h$ matrix and the $N^h$-vector with elements
\begin{equation}\label{eq:globalmatrix}
\left\{\begin{aligned}
(\mathbb{A}_{single})_{ij} = \mathcal{A}(\phi_j,\phi_i)
\qquad\qquad\qquad\qquad\,\,\,\,
&\qquad\mbox{for}\,\,\,i,j=1,\ldots,N^h
\\
(\vecbf{b}_{single})_i = \mcF(\phi_i)
-
\sum_{j=N^h+1}^{\widetilde N^h } \mathcal{A}(\phi_j,\phi_i) \vecgp_j
&\qquad
\mbox{for}\,\,
i=1,\ldots,N^h,
\end{aligned}\right.
\end{equation}
respectively. Then, the discrete FE problem \eqref{eq:weak_nonlochh} is equivalent to the linear algebraic system
\begin{equation}\label{eq:global_sys}
\mathbb{A}_{single} \vecu = \vecb_{single}
\end{equation}
for the unknown nodal coefficient vector $\vecbf{u}$. The matrix $\mathbb{A}_{single}$ is symmetric and, owing to the fact that {$\Gamma$ has nonzero measure in $\Re^2$,} it is also positive definite \cite{aksoylu2011variational,acta20,du19,du12,du13}.

For $v^h\in W^h$, we define the discrete energy functional  
\begin{align}\label{eq:disc_en_global}
    \mathcal{E}_{single}^h(v^h) = \frac{1}{2}\mathcal{A}(v^h,v^h) - \mcF(v^h). 
\end{align}
Then, the discrete nonlocal volume-constrained problem \eqref{eq:weak_nonlochh}, respectively \eqref{eq:global_sys}, can be expressed in terms of the equivalent minimization problem \cite{acta20,du19,du12,du13}
\begin{equation}
\begin{aligned}
&\mbox{\em the vector $\vecbf{u}$ solves \eqref{eq:global_sys}} 
\iff 
\mbox{\em $u^h(\xb)\in W^h$ solves \eqref{eq:weak_nonlochh}}
\iff
\\&\qquad 
\mathcal{E}^h_{single}(u^h(\xb)) = \min\limits_{v^h(\xb)\in W^h} {\mathcal E}^h_{single}( v^h(\xb)) \quad\mbox{\em subject to}\,\,\, v^h(\xb)=g^h(\xb) \mbox{ on $\Gamma$}.
\end{aligned}
\end{equation}

\begin{remark}\label{sissue}
{\em
As implied by \eqref{eq:int_dom}, any point $\xb\in\Omega\cup\Gamma$ interacts only with points in the ball $B_\delta(\xb)$. This raises a serious issue in FE methods for nonlocal problems because the intersection of {these} balls {with} the finite element grid {produces} cut elements, i.e., partial elements, within the ball. As a result, one either has to deal with cut elements or, if one wants to only deal with uncut elements, one is faced with discontinuous integrands that vanish outside the ball. This issue is glossed over in many FE papers for nonlocal problems, especially those that only provide one-dimensional numerical results. However, we address this issue in Section \ref{defsubgrid}. A comprehensive discussion of how to effectively handle cut elements can be found in \cite{DEliaFEM2020}.
}
\qquad$\Box$
\end{remark}

%%%%%%%%%%%%%%

\section{Nonlocal domain decomposition in the continuous setting}\label{sec:continuousDD}

In this section we first describe, in the continuous setting, how to define a nonlocal decomposition of the domain and then introduce the formulation of the multi-domain system. The central result of this section {proves the equivalence of the solution of the single-domain system and the one corresponding to the multi-domain system}. The significance of this result is that it establishes the consistency and the well posedness of our multi-domain formulation.

%%%%%%%%%%%%
\subsection{Construction of the geometric domain decomposition}\label{sec222111}

\begin{figure}[h!]
   \centering
   \includegraphics[width=2.25in]{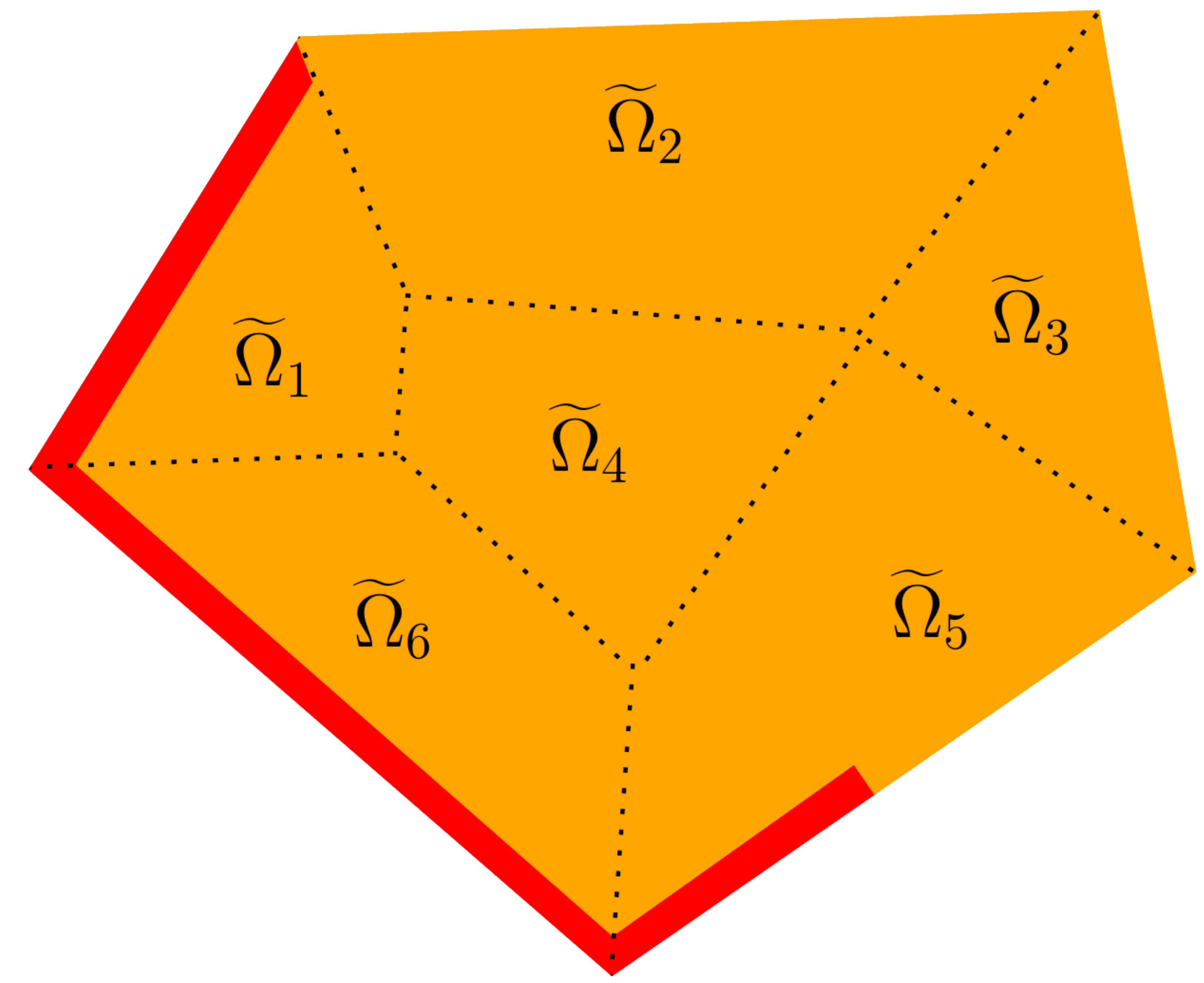}
   \quad 
   \includegraphics[width=2.25in]{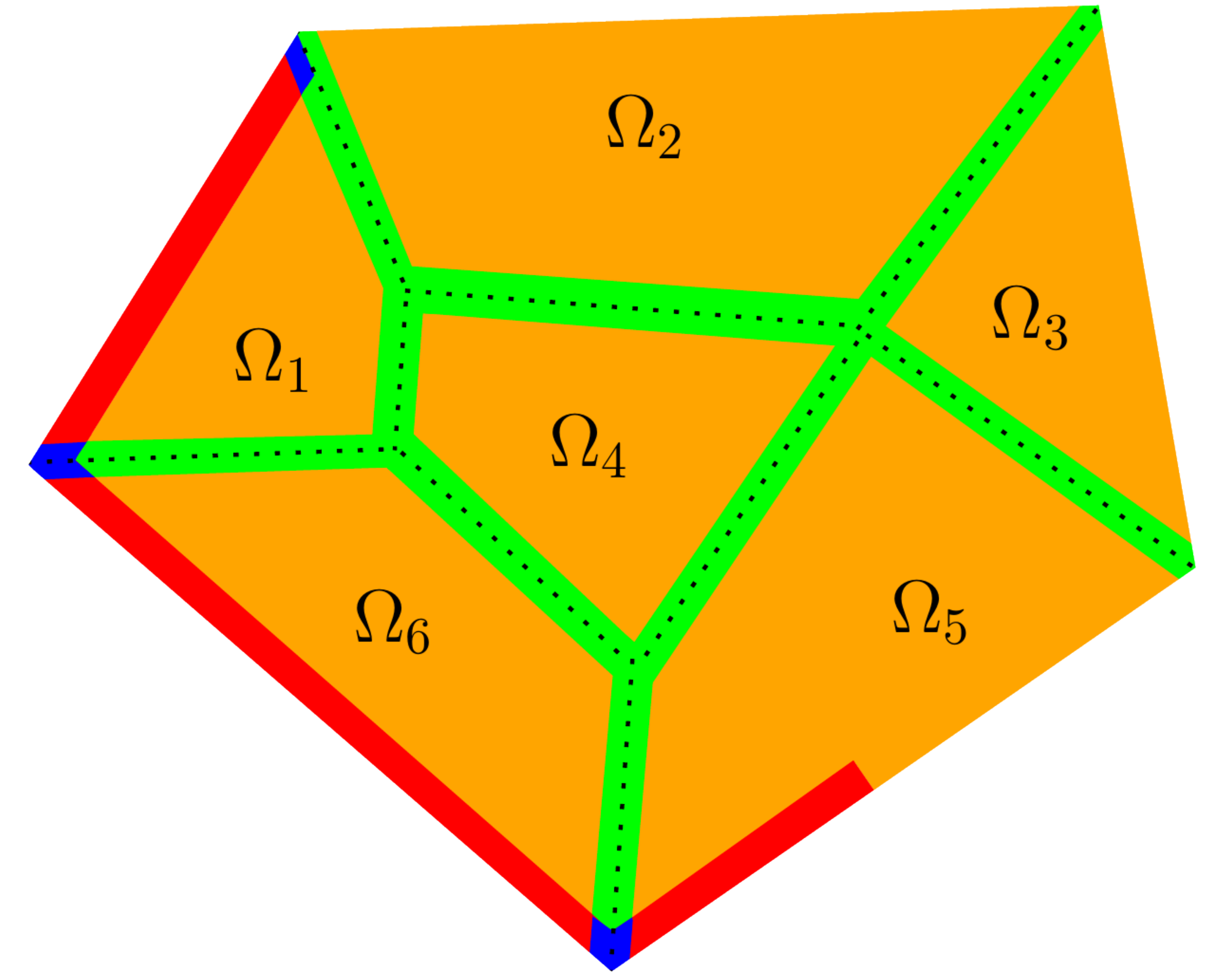}
\caption{Left: a non-overlapping, covering subdivision of the domain $\Omega$ into six subdomains $\widetilde\Omega_n$, $n=1,\ldots,N_s=6$. Right: the corresponding nonlocal overlapping domain subdivision $\Omega_n\cup\widehat\Gamma_n$, $n=1,\ldots,N_s=6$. The domain $\cup_{n=1}^{6}\widehat\Gamma_n$ is depicted in green and the domain $\Gamma=\Gamma_1\cup\Gamma_2\cup\Gamma_5\cup\Gamma_6$ is depicted in red with $\Gamma_3$ and $\Gamma_4$ being empty sets. The blue regions in $\Gamma$ illustrate the overlaps between pairs of $\Gamma_n$.}
   \label{fig2.2}
\end{figure}

In a standard PDE domain decomposition setting, one can partition $\Omega$ into non-overlapping subdomains and then simply define the subdomain problems by restricting the global operator to each subdomain. Such a construction is impossible in the nonlocal setting due to the inherent nonlocal interactions which require any two adjacent subdomains to share an interface having nonzero volume. As a result, our substructuring-based domain decomposition starts from a non-overlapping, covering subdivision of $\Omega$ into $N_s$ subdomains $\{\widetilde\Omega_n\}_{n=1}^{N_s}$, as illustrated in Figure \ref{fig2.2}-left for $N_s=6$, and then adds the overlaps necessary for the nonlocal interactions. 
Note that some of the domains $\widetilde\Omega_n$ include part of the boundary $\partial\Omega$. For example, in Figure \ref{fig2.2}-left, we have that this is the case for $\widetilde\Omega_2$, $\widetilde\Omega_3$, and $\widetilde\Omega_5$.
%%%

For each subdomain $\widetilde\Omega_n$, $n=1,\ldots,N_s$, we define the (possibly semi-) open smaller subdomain 
\begin{equation}\label{subdot}
   \Omega_n =
   \big\{ \xb\in\widetilde\Omega_n\,\,\,
   :\,\,\, |\yb-\xb|>\frac\delta2 \quad\forall \,\yb\in
 \Omega\setminus\widetilde\Omega_n  
\big\};\end{equation}
see Figure \ref{fig2.2}-right for an illustration. Note that $\{\Omega_n\}_{n=1}^{N_s}$ is, by construction, a set of {\em non-overlapping} domains. We also subdivide the interaction domain $\Gamma$ into a set of {\em overlapping}, covering subdomains 
{\begin{equation}\label{subdg}
    \Gamma_n = \big\{ \xb\in\Gamma\,
    : \, \exists\;\yb\in \Omega_n \;\; \hbox{such that}
    \;\;|\yb-\xb| \le \delta 
    \;\big\}\qquad\mbox{for $n=1,\ldots,N_s$}
\end{equation}}
and also define the set of {\em overlapping} subdomains
{\begin{equation}\label{subdgt}
    \widehat\Gamma_n = \big\{ \xb\in \Omega\setminus\Omega_n\,
    : \, \exists\;\yb\in \Omega_n \;\; \hbox{such that}
    \;\;|\yb-\xb| \le \delta 
    \;\big\}\qquad\mbox{for $n=1,\ldots,N_s$}.
\end{equation}}

In Figure \ref{fig2.2}-right, for each $n$, $\Gamma_n\cup \widehat\Gamma_n$ consists of all the strips of thickness $\delta$ that surround $\Omega_n$, including in some instances a portion of $\Gamma$. Because each subdomain $\widehat\Gamma_n$ overlaps with at least one other subdomain $\widehat\Gamma_{n'}$, $n'\ne n$, {and some subdomain $\Gamma_n$ may overlap with another subdomain $\Gamma_{n'}$, $n'\neq n$}, the above construction results in the {\em overlapping} domain decomposition of $\Omega\cup\Gamma$ given by
\begin{equation}\label{overdd}
\Omega\cup\Gamma = \cup_{n=1}^{N_s} \Omega_n\cup\widehat\Gamma_n\cup\Gamma_n.
\end{equation}

Following conventional DD nomenclature, we subdivide the set of subdomains $\{\Omega_n\cup\widehat\Gamma_n\cup\Gamma\}_{n=1}^{N_s}$ into two classes:
$$
\begin{aligned}
\mbox{\em floating subdomains} & \quad\mbox{if \,\,\,$\Gamma_n=\emptyset$}
\\
\mbox{\em non-floating subdomains} & \quad\mbox{if \,\,\,$\Gamma_n \ne \emptyset$.}
\end{aligned}
$$ 
For example, in Figure \ref{fig2.2}-right, $\Omega_3$ and $\Omega_4$ are floating subdomains, whereas $\Omega_1$, $\Omega_2$, $\Omega_5$, and $\Omega_6$ are non-floating. 
Analogous to the conventional local DD setting, a floating domain is endowed with a purely Neumann nonlocal volume constraint so that its associated nonlocal problem has a non-trivial null space. In a typical local or nonlocal DD configuration, the number of subdomains is large and most of them are of the floating type.

\begin{remark}\label{overlapdom}
{\em As alluded to in Section \ref{intro}, in a local PDE setting one can consider both non-overlapping and overlapping DD algorithms because the latter offer some computational conveniences and may allow for a faster convergence of iterative solution methods. Typically, the overlap depends on the grid size and its thickness goes to zero as the mesh size is reduced. In contrast, the nonlocal setting \emph{requires} any two adjacent subdomains to overlap in order to compute the necessary nonlocal interactions between them. As a result, the size of this \emph{mandatory} overlap is determined not by the mesh size but by the interaction radius $\delta$, i.e., its thickness is independent of the underlying discretization mesh.
}\qquad$\Box$
\end{remark}

%%%%%%%%%%%%
\subsection{The domain decomposition (multi-domain) system}

A multi-domain system is a system of $N_s$ equations, each of which holds for $\xb\in\Omega_n\cup\widehat\Gamma_n\cup\Gamma_n$, $n=1,\ldots,N_s$. In constructing the multi-domain system we have to deal with overlapping domains, i.e., although the three domains $\cup_{n=1}^{N_s}\Omega_n$, $\cup_{n=1}^{N_s}\widehat\Gamma_n$, and $\cup_{n=1}^{N_s}\Gamma_n$ are mutually disjoint as are the $N_s$ domains $\Omega_n$, {there are overlaps among the $N_s$ domains in $\cup_{n=1}^{N_s}\widehat\Gamma_n$ and there could be overlaps among the $N_s$ domains in $\cup_{n=1}^{N_s}\Gamma_n$.}

%%%%%%%%%%%%%%%%%%%%%%%%%%%%%%%%%%%%%%%%%%%%%%%%%%
%\subsubsection{Dealing with overlapping subdomains}\label{domoverl}

To properly deal with the consequences of having overlapping domains, we define the {generalized characteristic}  functions
\begin{equation}\label{eq:zeta}
\zeta_{\mcA}(\xb,\yb)=
\sum_{n=1}^{N_s}\mcX_{\Omega_n\cup\widehat\Gamma_n\cup\Gamma_n}(\xb)
\mcX_{\Omega_n\cup\widehat\Gamma_n\cup\Gamma_n}(\yb)   
\quad
\mbox{and}
\quad
\zeta_{\mcF}(\xb)= \displaystyle\sum_{n=1}^{N_s}\mcX_{\Omega_n\cup\widehat\Gamma_n}(\xb).
\end{equation}
Note that $\zeta_{\mcA}(\xb,\yb)$ is a symmetric function, i.e., $\zeta_{\mcA}(\xb,\yb)=\zeta_{\mcA}(\yb,\xb)$ and is a non-negative piecewise integer-valued function and  $\zeta_{\mcF}(\xb)$ is a positive piecewise integer-valued function. 
\begin{remark}\label{zetazero}
{\em It is possible for $\zeta_{\mcA}(\xb,\yb)=0$. For example, this is the case if $\xb\in\Omega_n$ and  $\yb\in\Omega_{n'}$ with $n'\ne n$. However, this {does not present a problem for} \eqref{eq:bic} below, {which involves $\zeta_{\mcA}(\xb,\yb)^{-1}$}, because points $\xb$ interact only with points $\yb$ such that $|\yb-\xb|\le \delta$ and for such pairs of points, $\zeta_{\mcA}(\xb,\yb)>0$.
}
\qquad$\Box$
\end{remark}

%%%%%%%%%%%%%%%%%%%%%%%%%%%%%%%%%%%%%%%%%%%%%%%%%%%
\subsubsection{The subdomain system}\label{subdsyst}

For $n=1,\ldots,N_s$ and any pair of functions $u_{n}(\xb)$ and $v_{n}(\xb)$ defined on $\Omega_n \cup \widehat{\Gamma}_n\cup\Gamma_n$, 
we define the subdomain bilinear form
\begin{equation}\label{eq:bic} 
\begin{aligned}
&\mathcal{A}_n(u_n,v_n)
\\&
=\int\limits_{\Omega_n\cup\widehat\Gamma_n\cup\Gamma_n}^{}\,\,    
 \int\limits_{\Omega_n\cup\widehat\Gamma_n\cup\Gamma_n}^{}
\hspace{-.2in} 
 \zeta_{\mcA}(\xb,\yb)^{-1} \big(u_n(\yb) - u_n(\xb)\big)\big(v_n(\yb) - v_n(\xb)\big)\gamma(\xb,\yb)d\yb  d\xb 
\end{aligned}
\end{equation}
and the associated subdomain linear functional
\begin{equation}\label{eq:lfc} 
\mcF_n(v_n) = \int_{\Omega_n\cup\widehat\Gamma_n} \zeta_{\mcF}(\xb)^{-1}v_n(\xb)f(\xb)d\xb,
\end{equation}
where, of course, for floating domains, i.e., if $\Gamma_n=\emptyset$, the integrals over $\Gamma_n$ vanish.

For $n=1,\ldots,N_s$, we define the function spaces
\begin{equation}\label{spacesdd}
\begin{aligned}
&\left\{
\begin{aligned}
&W_n = \{w\in L^2 (\Omega_n\cup\widehat\Gamma_n\cup\Gamma_n) \,\,:\,\, |||w|||_n < \infty\}
\\&
{W_n^0} = \{w\in W_n\,\,:\,\, w=0 \,\,\mbox{for $\xb\in\Gamma_n$}\},
\end{aligned}
\right.
\\& 
\qquad\qquad\qquad\qquad
\mbox{where} \quad |||w|||_n^2 =  \mathcal{A}_n(w,w) + \|w\|^2_{L^2(\Omega_n\cup\Gamma_n\cup\widehat\Gamma_n)}.
\end{aligned}
\end{equation}
Let
\begin{equation}\label{wwtt}
\widetilde W_n = \left\{
\begin{aligned}
&W_n^0 \quad\mbox{if $\Gamma_n\ne\emptyset$}, \,\,\,\,\mbox{i.e., for non-floating domains}
\\
&W_n \quad\mbox{if $\Gamma_n=\emptyset$},\,\,\,\,\mbox{i.e., for floating domains.}
\end{aligned}
\right.
\end{equation}
We then define the {\em domain-decomposition} or {\em multi-domain} system of equations as
\begin{equation}\label{ddsystems}
\left\{\begin{aligned}
&\mbox{\em given $f(\xb)\in W'$, $g(\xb)\in W_\Gamma$, and a kernel $\gamma(\xb,\yb)$,}
\\&
\mbox{\em for $n=1,\ldots, N_s$, }\mbox{\em find $u_n\in W_n$ such that}
\\&
\quad
\begin{aligned}
&\mathcal{A}_n(u_n,v_n) = \mcF_n(v_n) 
\quad
\forall\, v_n\in \widetilde W_n
\end{aligned}
\hspace{2.05in}(a)
\\
&\mbox{\em subject to}
\\
&\quad
 u_n(\xb) = u_{n'}(\xb) \,\,
     \quad \mbox{$\forall\,\xb\in \widehat\Gamma_{n} \cap \widehat{\Gamma}_{n'}$\,\,  \em for $n'=n+1,\ldots,N_s$}   
\hspace{.8in}(b)
\\
&\mbox{\em and}
\\
&\quad u_n(\xb) = g(\xb) \,\,
     \quad \forall\,\xb\in\Gamma_{n}
     \quad\mbox{\em if $\Gamma_n\ne\emptyset$}, \,\,\mbox{\em i.e., for non-floating domains}.
     \hspace{.1in}(c)
\end{aligned}\right.
\end{equation}

Equation (\ref{ddsystems}b) can be thought of as a nonlocal version of the standard continuity constraint in non-overlapping local DD methods. 
The constraints in that equation are needed
because the solution $u(\xb)$ of \eqref{eq:weak_nonloc} is a single-valued function on $\Omega\cup\Gamma$, and in particular on $\cup_{n=1}^{N_s} \widehat\Gamma_n$. However, by construction, we have that for all $n'$ such that $\widehat\Gamma_n\cap\widehat\Gamma_{n'}\ne\emptyset$, both $u_n(\xb)$ and $u_{n'}(\xb)$ are defined on $\widehat\Gamma_n\cap\widehat\Gamma_{n'}\ne\emptyset$. Clearly, we have that (\ref{ddsystems}b) must be imposed on that domain. 
Of course, this equation automatically holds on $\Gamma_n\cap\Gamma_{n'}$ because both $u_n(\xb)=g(\xb)$ and $u_{n'}(\xb)=g(\xb)$ there.
\begin{remark}\label{weakconst}
{\em As in the standard (local) DD case, the constraints in (\ref{ddsystems}b) ensure single-valued solutions of the multi-domain system, and are appropriate  when the global nonlocal solution is also continuous. However, one of the principal advantages of nonlocal models is that, for some kernels in common use, they admit solutions with jump discontinuities  \cite{acta20,du19,du12,du13}. {To allow for} such solutions, one may choose to enforce (\ref{ddsystems}b) weakly, i.e., for $n=1,\ldots,N_s$ and $n'=n+1,\ldots,N_s$,
\begin{equation}\label{weakconsteq}
\int_{\widehat\Gamma_{n} \cap \widehat{\Gamma}_{n'}}
\big( u_n(\xb) - u_{n'}(\xb) \big)v(\xb) d\xb = 0
\qquad \forall\, v\in W|_{\widehat\Gamma_{n} \cap \widehat{\Gamma}_{n'}},
\end{equation}
Note that \eqref{weakconsteq} also arises when (\ref{ddsystems}b) is enforced using Lagrange multipliers; such a treatment of (\ref{ddsystems}b) would {result in finite element formulations of the constraints (\ref{ddsystems}b) involving mass matrices}.
}
\qquad$\Box$
\end{remark}

\begin{remark}\label{redunconst}
{\em The constraints in (\ref{ddsystems}b) are not independent. For example, consider a point 
$\xb\in\widehat\Gamma_1\cap\widehat\Gamma_4\cap\widehat\Gamma_6$ near the bottom right corner  of $\Omega_1$ in Figure \ref{fig2.2}-right. Then, (\ref{ddsystems}b) would include the constraints $u_1(\xb)= u_4(\xb)$, $u_1(\xb)= u_6(\xb)$, and $u_4(\xb)= u_6(\xb)$, only two of which are independent. These redundancies in (\ref{ddsystems}b) have implications in the design of discretization algorithms as is discussed in Section \ref{multifem}.
}
\qquad$\Box$
\end{remark}

%%%%%%%%%%%%%%%%%%

\subsubsection{Equivalence of the single-domain and multi-domain problems}\label{equivsdmd}

Our next task is to show that the solution $u(\xb)$ for $\xb\in\Omega\cup\Gamma$ of the single-domain system \eqref{eq:weak_nonloc} and the solutions $u_n(\xb)$, $n=1,\ldots,N_s$, of the multi-domain system \eqref{ddsystems} are the same or, more precisely, that $u_n(\xb)=u(\xb)$ for $\xb\in\Omega_n\cup\widehat\Gamma_n\cup\Gamma_n$. 
The first step towards that end is the following lemma. We refer to Appendix \ref{main-proof} for a proof.
%%%
\begin{lemma}\label{lemma111}
Given functions $u(\xb)$ and $v(\xb)$ for $\xb\in\Omega\cup\Gamma$, define the functions $u_n(\xb)$ and $v_n(\xb)$ for $\xb\in \Omega_n\cup\widehat\Gamma_n\cup\Gamma_n$ as
\begin{equation}\label{uequn}
  u_n(\xb) = u(\xb)\big|_{\Omega_n\cup\widehat\Gamma_n\cup\Gamma_n}
  \quad\mbox{and}\quad
  v_n(\xb) = v(\xb)\big|_{\Omega_n\cup\widehat\Gamma_n\cup\Gamma_n} 
  \quad\mbox{for $n=1,\ldots,N_s$}.
\end{equation}
Then, for $n=1,\ldots,N_s$ and $n'=n+1,\ldots,N_s$, we have that, if $\widehat{\Gamma}_n\cap \widehat{\Gamma}_{n'}\ne\emptyset$ and ${\Gamma}_n\cap {\Gamma}_{n'}\ne\emptyset$,
\begin{equation}\label{unequnp}
     u_{n'}(\xb) = u_{n}(\xb) 
     \,\,\,\,\mbox{and}\,\,\,\, v_{n'}(\xb) = v_{n}(\xb) 
     \quad \mbox{for} \,\,\,\, \xb\in \widehat{\Gamma}_n\cap \widehat{\Gamma}_{n'}
     \,\,\,\,\mbox{and}\,\,\,\,
     \xb\in {\Gamma}_n\cap {\Gamma}_{n'}.
\end{equation}
Then, for the bilinear forms and linear functionals defined in \eqref{eq:bil1}, \eqref{eq:bic}, and \eqref{eq:lfc}, we have that
\begin{equation}\label{suman}
    \mathcal{A}(u,v) =\displaystyle\sum_{i=1}^{N_s} \mathcal{A}_n(u_n,v_n) 
\qquad\mbox{and}\qquad
 \mcF(v) =  \displaystyle\sum_{n=1}^{N_s} \mcF_n(v_n).
\end{equation}
\end{lemma}
%%%

Using Lemma \ref{lemma111} and the fact that $W = \otimes_{n=1}^{N_s} W_n$ and $W^0 = \otimes_{n=1}^{N_s} \widetilde W_n$ one can easily prove the following equivalence result.

%%%
\begin{proposition}\label{prop-equivalent}
Let $u(\xb)\in W$ denote the solution of the single-domain system \eqref{eq:weak_nonloc} and, for $n=1,\ldots,N_s$, let $u_n(\xb)\in W_n$ denote the solution of the $n$-th subproblem in multi-domain system \eqref{ddsystems}. Then, the multi-domain system \eqref{ddsystems} and the single-domain system \eqref{eq:weak_nonloc} are equivalent and their respective solutions coincide, i.e, $u_n(\xb) = u(\xb)\big|_{\Omega_n\cup\widehat\Gamma_n\cup\Gamma_n}$.  
\hfill$\Box$
\end{proposition}

The equivalence of the multi-domain and the global problems and the fact that the latter is well posed imply that \eqref{ddsystems} is also well posed.

We note that the equivalence between the multi-domain weak formulation \eqref{ddsystems} and the single-domain system \eqref{eq:weak_nonloc} can also be established by noting that the single-domain energy functional \eqref{singleenergy} can be written as a sum of subdomain energy functionals defining a multi-domain energy functional. The proof is straightforward and is omitted.

\begin{proposition}\label{prop2.1}
Define the subdomain energy functionals $\mathcal{E}_n(u_n)$ by
\begin{equation}\label{ddenergies}
  \mathcal{E}_n(u_n):= \frac{1}{2}\mathcal{A}_n(u_n,u_n) - \mcF_n(u_n) \qquad  \mbox{for $n=1,\ldots,N_s$},
\end{equation}
where the bilinear form ${\mcA}_n(\cdot,\cdot)$ and linear functional $\mcF_n(\cdot)$ are defined in \eqref{eq:bic} and \eqref{eq:lfc}, respectively. Then,
\begin{align}\label{eq:functional_sum}
\displaystyle\sum_{n=1}^{N_s}\mathcal{E}_n(u_n) = \mathcal{E}_{single}(u).
\end{align}
Furthermore, the multi-domain weak formulation \eqref{ddsystems} is the Euler-Lagrange equation corresponding to the minimization problem
\begin{equation}\label{eq:4energies}
\inf_{v_n\in W_n,\, n=1,\ldots,N_s}
   \sum\limits_{n=1}^{N_s} \mathcal{E}_n(v_n)
   \quad
\left\{\begin{aligned}
&\mbox{subject to, for $n,n'=1,\ldots,N_s$, $n'\ne n$,}
\\&\mbox{$u_n(\xb) = u_{n'}(\xb)$ for $\xb\in \widehat\Gamma_{n} \cap \widehat{\Gamma}_{n'}$ and}
 \\&   \mbox{$u_n(\xb) = g(\xb)$ for $\xb\in \Gamma_n$. }  
\end{aligned}\right.
 \end{equation}  
\end{proposition}

%%%%%%%%%%%%
\section{Finite element discretization of the subdomains systems}\label{femsubsys}

In this section we {describe} a finite element discretization {of} the multi-domain system \eqref{ddsystems}. {We follow a standard FE assembly procedure specialized to our needs.} 
{Specifically,  for each $n=1,\ldots,N_s$, we proceed as follows.}

\hangone{\em { 1. Subdomain meshing.}}~{This step introduces a suitable finite element partition on}  each subdomain $\Omega_n\cup\Gamma_n\cup\widehat\Gamma_n$.

\hangone{\em {2. Selection of a subdomain finite element space.}}
{We endow each subdomain $\Omega_n\cup\Gamma_n\cup\widehat\Gamma_n$  with a finite element space $\widetilde W_n^h$ defined as the span of a set of basis functions}, usually chosen to be piecewise polynomials with respect to {the underlying grid.}  

\hangone{\em {3. Subdomain finite element ansatz}.} {We seek an approximation of the solution $u_n(\xb)$ of \eqref{ddsystems} out of the finite element space $\widetilde W_n^h$.}

\hangone{\em {4. Subdomain finite element assembly.}} {We construct the stiffness matrix and right-hand side vector for each subdomain in the same manner as for the global problem, i.e., by inserting the finite element basis functions into the respective weak subdomain equations.}

{As is the case with any DD method, the goal is to perform these steps in such a way so as to ensure that}
\begin{equation}\label{goalgoal}
\begin{aligned}
&\mbox{\em the {global solution obtained from} the $N_s$  {\em discretized} subdomain FE}
\\[-.75ex]
&\mbox{\em  systems corresponding to \eqref{ddsystems}, i.e., $u^{dd,h}(\xb)$, should be the same as}
\\[-.75ex]
&\mbox{\em the solution $u^h$ of the discretized single-domain FE system \eqref{eq:weak_nonlochh}.}
\end{aligned}
\end{equation}

{A standard way to achieve this goal in DD methods is to use a partition based on the underlying mesh. In this approach the subdomains are defined as unions of adjacent elements using a graph partitioning tool. The associated finite element spaces are then simply restrictions of the global space to the resulting subdomain. The first step, i.e., the subdomain meshing in this case is trivial because it uses an already defined global mesh. 
This construction is not restrictive as we are interested in DD as means to devise fast parallel solution algorithms for a given global nonlocal problem. In such a case a global mesh and a finite element space can be assumed to have been already constructed. This should be contrasted with mortar methods and mesh tying where from the onset one assumes that each subdomain is meshed separately.}

%%%%%%%%%%%%%%
\subsection{Definition of subdomain grids}\label{defsubgrid}

{We assume the setting of Section \ref{sec1-1}, i.e., that we are given a finite element mesh $\mcT^h$ of $\Omega\cup\Gamma$ that conforms to the interface shared by $\Omega$ and $\Gamma$. As before, $\widetilde N^h$ denotes the total number of nodal degrees of freedom on this mesh and $N^h$ are the nodal values associated with $\Omega$. The remaining $\widetilde N^h-N^h$ degrees of freedom live on the closed domain $\Gamma$.}

In {\em local} PDE settings, for both non-overlapping and overlapping DD methods, the next step is to subdivide the domain $\Omega$ into subdomains that contain only whole finite elements. If one is going to invoke a non-overlapping DD method, the construction of the subdomains is complete. If instead an overlapping DD method is to be used, one adds, to each non-overlapping subdomain, whole elements in neighboring subdomains that are within a certain distance from the common boundary between the two domains; the distance used is usually related to some multiple of the local grid size, although other criteria are also in use \cite{cai,mathew2008domain,toselli2006domain}. 
 
In the {\em nonlocal} case, the practical construction of a subdivision of $\Omega$ is similar to that for overlapping DD in the local case. We again start by subdividing $\Omega$ into the set $\{\widetilde\Omega_n\}_{n=1}^{N_s}$ of subdomains with each $\widetilde\Omega_n$ consisting of whole finite elements. We now want to add to and subtract from each subdomain $\widetilde\Omega_n$ strips of thickness $\delta/2$ to create the subdomains $\Omega_n$ and $\widehat\Gamma_n$. Unfortunately, because we are given a global grid over $\Omega\cup\Gamma$ to work with, in general, we will not be able to (see Remark \ref{sissue}) define grids consisting of whole finite elements that respect the boundaries between $\Omega_n$ and $\widehat\Gamma_n$, i.e., those domains would also contain partial (cut) elements which is something we want to avoid. 

We thus see that there is a big difference between the definitions of overlapping grids for the local and nonlocal cases. To recapitulate, in the local case we are free to add whole elements to effect an overlap with neighboring elements. In the nonlocal case, we do not have this freedom because the strips to be created have thickness $\delta$ irrespective of the given global finite element grid so that, in general, that strip will not consist of whole elements. Thus we have the choice of truncating {elements} so that the $\delta$ thickness of the strip is respected or instead {\em approximate} the strip by a strip consisting of whole {elements} which is tantamount to approximating the common boundaries of $\Omega_n$ and $\widehat\Gamma_n$ by element edges. We use the latter choice because it is substantially easier to implement and, as shown below, does not compromise achieving the goal \eqref{goalgoal}.

The above discussion motivates the following procedure for the construction, in the discretized setting, of a subdivision of $\Omega$ into subdomains that is analogous, but not the same, as that in Section \ref{sec222111} for the continuous problem. We begin by assuming, as is done in Section \ref{sec1-1} for the single-domain setting, that 

\hangone~{--}~we are given an integer $N_s>1$ and a finite element meshing $\mcT^h$ of $\Omega\cup\Gamma$ which respects their common boundary $\overline\Omega\cap\Gamma$. 

\noindent We denote by $\mcT^h_\Omega$ and $\mcT^h_\Gamma$ the sets of finite elements in $\Omega$ and $\Gamma$, respectively, and we denote by $\TTT$ a typical element in $\mcT^h$. Then, 

\hangone~{--}~we subdivide $\Omega$ into $N_s$ non-overlapping, covering subdomains $\widetilde\Omega_n$, $n=1,\ldots,N_s$, such that {\em each subdomain $\widetilde\Omega_n$ consists entirely of whole finite elements.}

\noindent This step is effected exactly in the same manner as for the local PDE non-overlapping DD setting so that no further comments are needed. 

Note that the boundary $\partial\widetilde\Omega_n$ of $\widetilde\Omega_n$ consists of two or three disjoint, covering parts. First, we have for all $n$, 

\hangtwo~Type 1.~$\bigcup_{n'=1,\,n'\ne n}^{N_s} \partial\widetilde\Omega_n\cap\partial\widetilde\Omega_{n'}$, i.e., the common boundary shared by $\widetilde\Omega_n$ and subdomains $\widetilde\Omega_{n'}$ that abut to $\widetilde\Omega_n$.

\noindent We also have either one or both of

\hangtwo~Type 2.~$\partial\widetilde\Omega_n\cap\partial\Gamma$, i.e., the common boundary shared by $\widetilde\Omega_n$ and $\Gamma$.

\hangtwo~Type 3.~$\partial\widetilde\Omega_n\setminus \left[(\bigcup_{n'=1,\,n'\ne n}^{N_s} \partial\widetilde\Omega_n\cap\partial\widetilde\Omega_{n'})\right] \cup(\partial\widetilde\Omega_n\cap\partial\Gamma)$, i.e., the part of $\partial\widetilde\Omega_n$ that is not shared with the boundary of $\Gamma$ or with any of the boundaries of other subdomains $\widetilde\Omega_{n'}$. 

\noindent For example, referring to Figure \ref{fig2.2}-left, we have that the boundaries of floating domains such as $\widetilde\Omega_{3}$ consist of only Type 1 and 3 parts, the boundaries of domains such as $\widetilde\Omega_{6}$ consist of only Type 1 and 2 parts, and the boundaries of domains such as $\widetilde\Omega_{5}$ consists of all three parts.

{Recall that, by construction, the domains $\widetilde\Omega_n$ consist of whole FE triangles. However, in general, the subdomains $\Omega_n$, $\widehat\Gamma_n$, and $\Gamma_n$ consist of whole FE triangles {\em and} additionally partial (cut) FE triangles. To obtain the equivalence of the single- and multi-domain FE solutions, we necessarily have to work with subdomains that consist of only whole FE triangles because only such triangles, i.e., triangles $T\in\mathcal T_h$, are used in the single-domain FE method. For this reason, we define an approach for the construction of subdomains $\Omega_n^h$, $\widehat\Gamma_n^h$, and $\Gamma_n^h$ of $\Omega\cup\Gamma$ in such a way that all subdomains consists of only whole FE triangles. Additionally, the construction process is required to account for all interactions that occur between two subdomains $\widetilde\Omega_n$ and $\widetilde\Omega_{n'}$. Meeting this requirement is guaranteed if all triangles that overlap with the $\widehat\Gamma_n$ are included in $\widehat\Gamma_n^h$.
}

{The specific geometric domain decomposition we use is defined as follows.
We denote by $\xb^{vertex}$ a typical vertex on $\partial\widetilde\Omega_n$ and by $\xb^{barycenter}_\TTT$ the barycenter of a typical finite element $\TTT\in\mcT^h$. Then, for $n=1,\ldots,N_s$, we define the subdomains
{\begin{equation}\label{ddsubdomainsdd}
\left\{\begin{aligned}
\widehat\Gamma_n^h &=
\big\{ \TTT\in \mcT^h_\Omega \,:\,
\exists\; \xb^{vertex} \in \mbox{Type 1 part of $\partial\widetilde\Omega_n$} 
\\
&\qquad\qquad\qquad\qquad\hbox{such that} \;\;
| \xb^{vertex} - \xb_\TTT^{barycenter}| \le \delta/2+h\big\} 
\\
\Omega_n^h &= \big\{  \TTT\in \widetilde\Omega_n\setminus 
(\widetilde\Omega_n\cap \widehat\Gamma_n^h )   \big\} 
\\[.5ex]
\Gamma_n^h & = \big\{ \TTT\in \mcT^h_\Gamma \,:\,
\exists \; \xb^{vertex} \in \mbox{Type 2 part of $\partial\widetilde\Omega_n$}\\
&\qquad\qquad\qquad\qquad\hbox{such that} \;\;
| \xb^{vertex} - \xb_\TTT^{barycenter} | \le \delta+h \big\}.
\end{aligned}\right.
\end{equation}}
We have that $\widehat\Gamma_n^h$ consists of all elements $\TTT\in \mcT^h_\Omega$ whose barycenters are within a distance $\delta/2+h$ of some element vertex on the Type 1 part of the boundary of $\widetilde\Omega_n$. Also, $\Gamma_n^h$ consists of all elements $\TTT\in \mcT^h_\Gamma$ whose barycenters are within a distance $\delta+h$ of some element vertex on the Type 2 part of the boundary of $\widetilde\Omega_n$. Note that this procedure guarantees that the true interface region $\widehat\Gamma_n$ is fully contained in the approximate interface region $\widehat\Gamma_n^h$.} 
%For obvious reasons, we refer to the approach we use as ``barycenter-based''.

\begin{figure}[h!]
\begin{center}
\begin{tabular}{cccccc}
   \includegraphics[width=.65in]{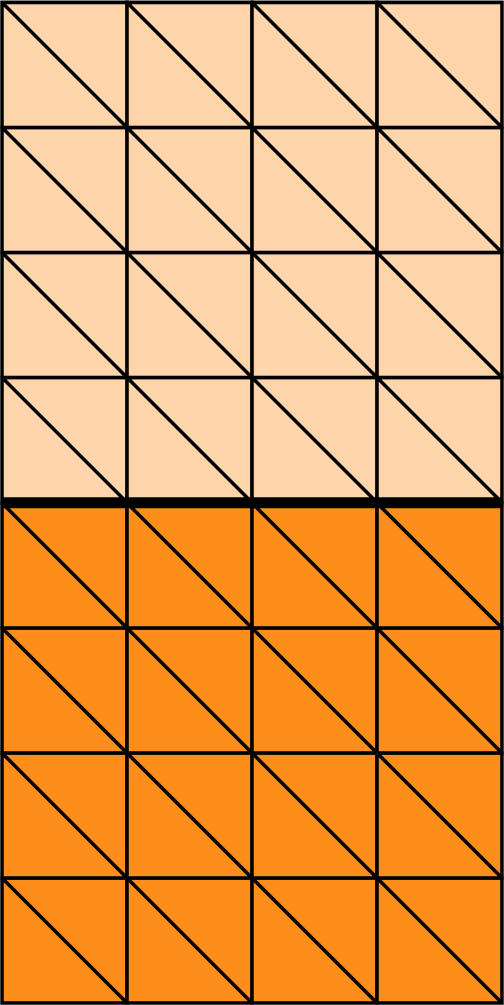}\,\,
&   \includegraphics[width=.65in]{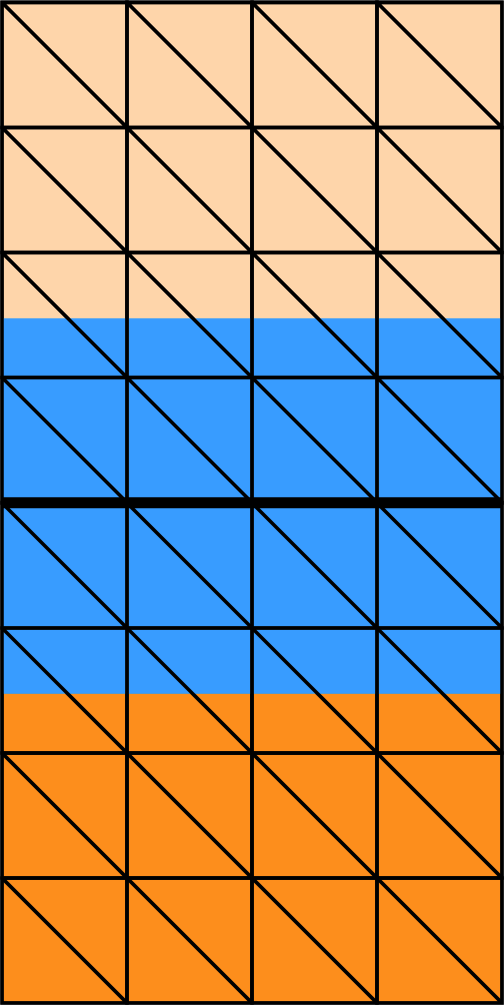}\,\,
 &   \includegraphics[width=.65in]{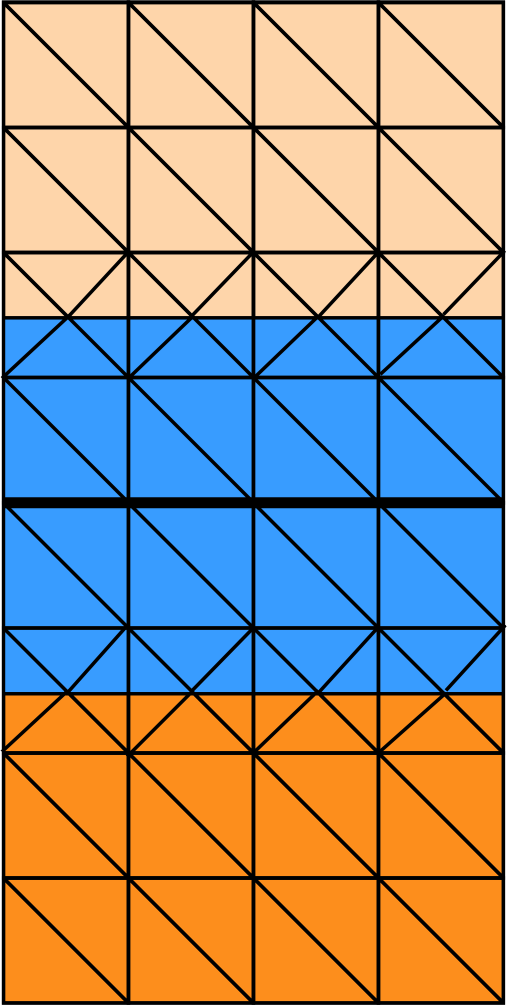}\,\,
  &   \includegraphics[width=.65in]{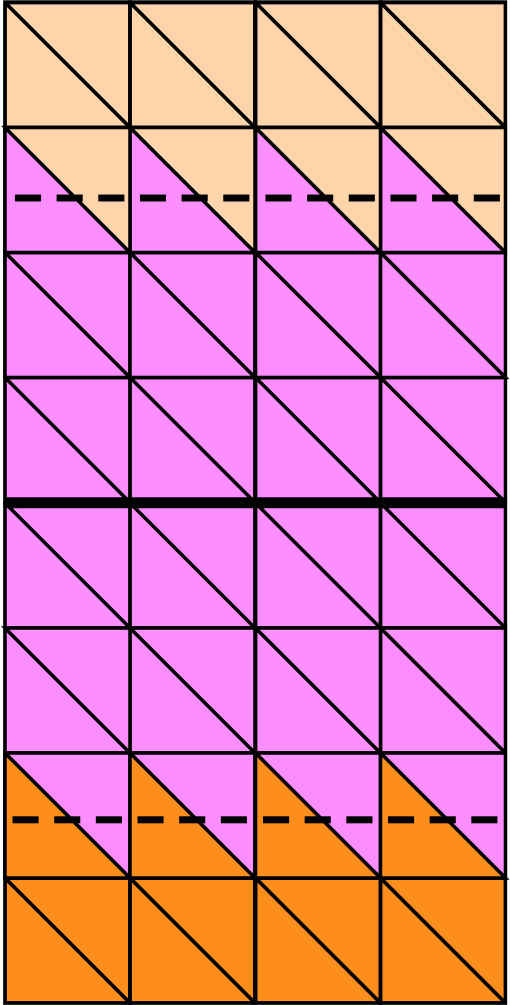}\,\,
 &   \includegraphics[width=.65in]{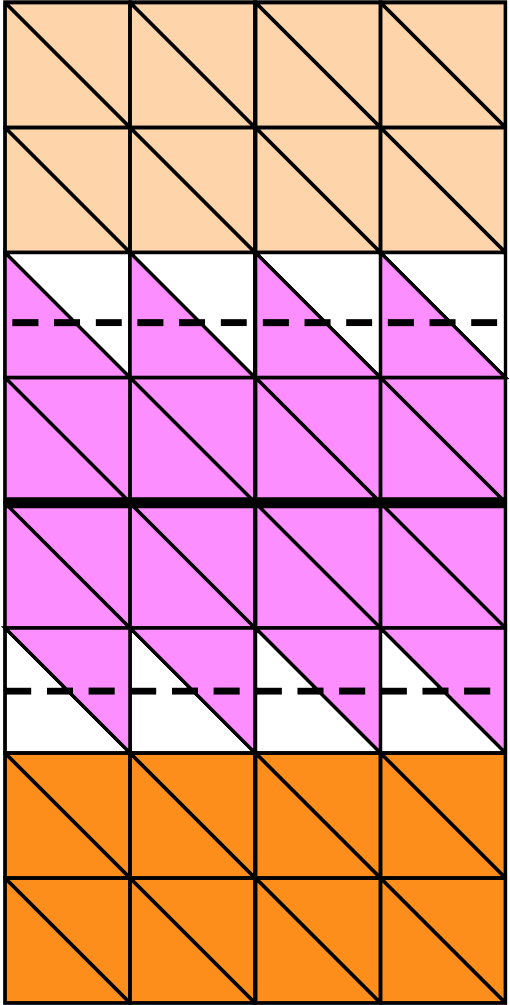}\,\,
 &   \includegraphics[width=.65in]{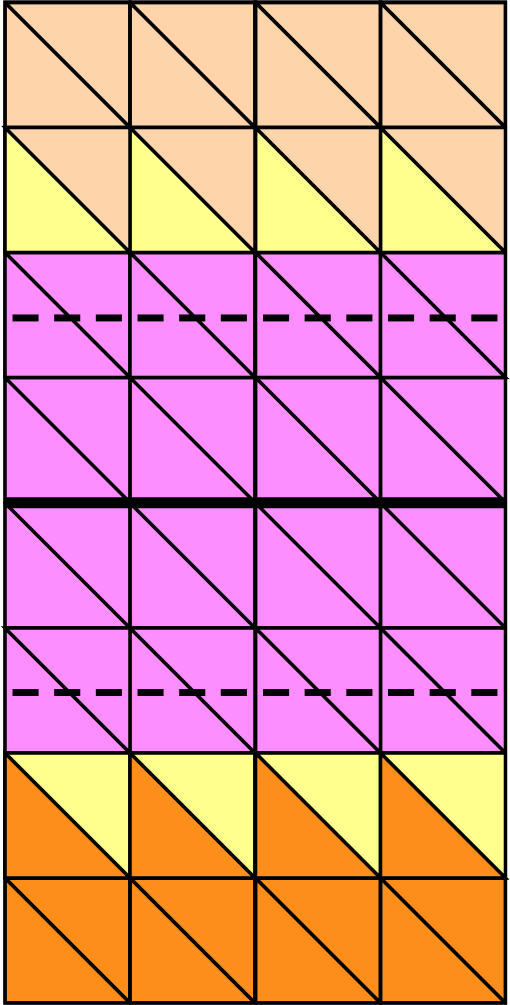}
\\
(a) & (b) & (c) & (d) & (e) & (f)
\end{tabular}
\end{center} 
\caption{(a): A portion of the FE single-domain grid and portions of the two subdomains $\widetilde\Omega_n$ and $\widetilde\Omega_{n'}$. (b): Portions of the subdomains {$\Omega_n$} and {$\Omega_{n'}$} (in shades of orange) and of $\widehat\Gamma_n\cap\widehat\Gamma_{n'}$ (in blue). (c): A re-meshing of the FE grid of (a) so that now the common boundaries between the subdomains in (b) are respected {but only if cut elements are introduced}. (d): Portions of the approximate subdomains {$\Omega_n^h$} and {$\Omega_{n'}^h$} (in shades of orange) and of $\widehat\Gamma_n^h\cap\widehat\Gamma_{n'}^h$ (in magenta) as determined using \eqref{ddsubdomainsdd}. {(e): The white elements overlap with $\widehat\Gamma_n$ {but are not included in the approximate interface generated with the criterion $| \xb^{vertex} - \xb_\TTT^{barycenter}| \le\delta/2$}. (f): The yellow elements are in $\widehat\Gamma_n^h\cap\widehat\Gamma_{n'}^h$ but do not overlap with $\widehat\Gamma_n\cap\widehat\Gamma_{n'}$.}
} 
   \label{dd-grid}
\end{figure}

We illustrate the above discussion in Figure \ref{dd-grid}. Note that all of that discussion applies even to the case of the single domain $\Omega$ being a rectangle and a single-domain FE grid that is Cartesian and uniform. Figure \ref{dd-grid}a depicts a portion of the  grid in the single domain $\Omega$ that respects the common boundary (depicted by the thick line segment) between two subdomains $\widetilde\Omega_n$ and $\widetilde\Omega_{n'}$. For Figure \ref{dd-grid}b, we have that the orange subdomains depict portions of the subdomains $\Omega_n$ and $\Omega_{n'}$ and the blue domain depicts a portion of $\widehat\Gamma_n\cap\widehat\Gamma_{n'}$. Note that the common boundaries of both $\Omega_n$ and $\Omega_{n'}$ with $\widehat\Gamma_n\cap\widehat\Gamma_{n'}$ do not respect the  grid so that $\Omega_n$, $\Omega_{n'}$, and $\widehat\Gamma_n\cap\widehat\Gamma_{n'}$ all contain some partial (cut) triangles. Figure \ref{dd-grid}c illustrates the need to make changes to the single-domain FE grid so that the new grid does respect those common boundaries. Of course, if we define the subdomain FE discretization using the new grid of Figure \ref{dd-grid}c, there is no hope for the solution of the FE discretization of \eqref{ddsystems} to be the same as the solution of single-domain FE discretization \eqref{eq:weak_nonlochh}, i.e., the goal \eqref{goalgoal} cannot be achieved. Note also that the re-meshing of Figure \ref{dd-grid} is relatively easy to effect for Cartesian grids, but becomes a much more complex task for general grids, especially in three dimensions. {Figure \ref{dd-grid}d illustrates the process defined in \eqref{ddsubdomainsdd}. Now the orange-shaded regions depict portions of the subdomains $\Omega_n^h$ and $\Omega_{n'}^h$ and the magenta region depicts a portion of $\widehat\Gamma_n^h\cap\widehat\Gamma_{n'}^h$. Note that, in Figure \ref{dd-grid}d, those three domains all contain only whole FE triangles. Also, the blue region in Figure \ref{dd-grid}b, i.e. $\widehat\Gamma_n\cap\widehat\Gamma_{n'}$, is fully contained in the set of magenta triangles.}
{Figure \ref{dd-grid}e illustrates why, e.g., in the first equation in \eqref{ddsubdomainsdd}, we used the criteria $| \xb^{vertex} - \xb_\TTT^{barycenter}| \le \delta/2+h$ and not simply the criteria $| \xb^{vertex} - \xb_\TTT^{barycenter}| \le \delta/2$. The barycenters of the white triangles are such that $| \xb^{vertex} - \xb_\TTT^{barycenter}| > \delta/2$ so that using the latter criteria means that the white triangles are not included in $\widehat\Gamma_n^h$ even though it is obvious from Figure \ref{dd-grid}b that those triangles overlap with $\widehat\Gamma_n^h\cap\widehat\Gamma_{n'}^h$. On the other hand, Figure \ref{dd-grid}f illustrates that using the criteria $| \xb^{vertex} - \xb_\TTT^{barycenter}| \le \delta/2+h$ results in the yellow triangles which are in $\widehat\Gamma_n^h\cap\widehat\Gamma_{n'}^h$ but do not overlap with $\widehat\Gamma_n\cap\widehat\Gamma_{n'}$ so that those elements do not interact with elements on the other side of $\widehat\Gamma_n^h\cap\widehat\Gamma_{n'}^h$. However, the equivalence between single- and multi-domain solutions is not compromised because in assembling the FE stiffness matrix, any additional entries in that matrix corresponding to points in the yellow elements are taken care of by weighting the integrand with the inverse of $\zeta_\mcA(\xb,\yb)$.}
%%%

\begin{remark}
{\em Comparing the definitions of $\widehat\Gamma_n^h$, $\Gamma_n^h$, and $\Omega_n^h$ with the definitions of $\widehat\Gamma_n$, $\Gamma_n$, and $\Omega_n$ given in Section \ref{sec222111}, one can certainly view the first trio as approximations to the second trio. However, this view does not intrude on any aspect of the developments that follow. For example, the accuracy of these domain approximations is, as is made evident in Section \ref{multifem}, irrelevant with respect to the goal stated in \eqref{goalgoal}.}\qquad$\Box$
\end{remark}

%%%%%%%%%%%%%%
\subsection{Multi-domain finite element system}\label{multifem}

Let $W^h$ denote the finite element space spanned by the set of basis functions $\{\phi_i(\xb)\}_{i=1}^{\widetilde N^h}$ used in Section \eqref{sec1-1} to define the single-domain FE system \eqref{eq:global_sys}. Then,  let  
\begin{equation}\label{wnhbasis}
W_n^h = \mbox{span} \big\{
%\cup_{i=1}^{\widetilde N^h}  
\phi_i(\xb) \,\,:\,\, {\xb_i}\in {\Omega_n^h\cup\Gamma_n^h\cup\widehat\Gamma_n^h} \big\}
\qquad\mbox{for $n=1,\ldots,N_s$},
\end{equation}
i.e., $W_n^h$ is spanned by {the basis functions $\phi_i$ such that the associated node $\xb_i$ belongs to  $\Omega_n^h\cup\Gamma_n^h\cup\widehat\Gamma_n^h$.}
The subspaces $W_n^h \subset W^h$ are well defined because of \eqref{ddsubdomainsdd}, i.e., because $\Omega_n^h$, $\Gamma_n^h$, and $\widehat\Gamma_n^h$ all consist of whole elements from the triangulation $\mcT^h$. We also define the spaces
\begin{equation}\label{wnhbasisc}
\begin{aligned}
{W_n^{0,h}} = \big\{w^h(\xb) &\in W_n^h\,\,\,:\,\,\,
w^h(\xb) = 0 \quad\mbox{for $\xb\in \Gamma_n^h$}\big\}
\\
&\qquad\mbox{for $n=1,\ldots,N_s$ such that $\Gamma_n^h\ne\emptyset$}.
\end{aligned}
\end{equation}
Note that {\em although we have that $ W_1^h\otimes\cdots\otimes  W_{N_s}^h = W^h\subset W$, in general, $ W_n^h\not\subset W_n$.}

We define
\begin{equation}\label{eq:zetah}
\zeta_{\mcA}^h(\xb,\yb)=
\sum_{n=1}^{N_s}\chi_{\Omega_n^h\cup\widehat\Gamma_n^h\cup\Gamma_n^h}(\xb)
\chi_{\Omega_n^h\cup\widehat\Gamma_n^h\cup\Gamma_n^h}(\yb)   
\quad
\mbox{and}
\quad
\zeta_{\mcF}^h(\xb)= \displaystyle\sum_{n=1}^{N_s}\chi_{\Omega_n^h\cup\widehat\Gamma_n^h}(\xb).
\end{equation}
For $n=1,\ldots,N_s$ and for $u_n^h(\xb),\,v_n^h(\xb)\in W_n^h $, we define the discretized subdomain bilinear form 
\begin{equation}\label{eq:bilf}
\begin{aligned}
  &\mathcal{A}_n^h(u_n^h,v_n^h) =
  \\& \int\limits_{\Omega_n^h\cup\Gamma_n^h\cup\widehat\Gamma_n^h}{}  \int\limits_{\Omega_n^h\cup\Gamma_n^h\cup\widehat\Gamma_n^h}{} 
\hspace{-.2in}  
  \zeta^h_{\mcA}(\xb,\yb)^{-1} \big(v_n^h(\yb) - v_n^h(\xb)\big)\big(u_n^h(\yb) - u_n^h(\xb)\big)\gamma(\xb,\yb)d\yb  d\xb
  \end{aligned}
  \end{equation}
and the discretized subdomain linear functional
\begin{equation}\label{eq:linf}
   \mcF_n^h(v_n^h) = \int_{\Omega_n^h\cup\widehat{\Gamma}_n^h} \zeta^h_{\mcF}(\xb)^{-1} v_n^h(\xb)f(\xb)d\xb 
\end{equation}
with the tacit understanding that for floating domains, i.e., if $\Gamma_n^h=\emptyset$, the domains of integration of both integrals in \eqref{eq:bilf} reduce to $\Omega_n^h \cup \widehat{\Gamma}_n^h$. 

{As done for the continuous problem, we introduce the function space}
{\begin{equation*}
\widetilde W^h_n = \left\{
\begin{aligned}
& W_n^{0,h} \quad\mbox{\em if $\Gamma_n^h\ne\emptyset$}, \,\,\,\,\mbox{\em i.e., for non-floating domains}
\\
&W_n^h \quad\mbox{\em\, if $\Gamma_n^h=\emptyset$},\,\,\,\,\mbox{\em i.e., for floating domains.}
\end{aligned}
\right.
\end{equation*}}
Then, we define the system of equations
\begin{equation}\label{ddsystemsh}
\left\{\begin{aligned}
&\mbox{\em given $f(\xb)\in W'$, $g(\xb)\in W_\Gamma$, and a kernel $\gamma(\xb,\yb)$,}
\\&
\mbox{\em for $n=1,\ldots, N_s$, }\mbox{\em find $u_n^h\in W_n^h$ such that}
\\&
\quad
\begin{aligned}
&\mathcal{A}_n^h(u_n^h,v_n^h) = \mcF_n^h(v_n^h) \qquad
\forall\, v_n^h\in  \widetilde W_n^h
\end{aligned}
\hspace{1.9in}(a)
\\
&\mbox{\em subject to}
\\
&\quad
 u_n^h(\xb) = u_{n'}^h(\xb) \,\,
     \quad \mbox{$\forall\,\xb\in \widehat\Gamma_{n}^h \cap \widehat{\Gamma}_{n'}^h$\,\,  \em for $n'=n+1,\ldots,N_s$}   
\hspace{.83in}(b)
\\
&\mbox{\em and}
\\
&\quad u_n^h(\xb) = g^h(\xb) \,\,
     \,\,\,\, \forall\,\xb\in\Gamma_{n}^h
     \quad\mbox{\em if $\Gamma_n^h\ne\emptyset$}, \,\,\mbox{\em i.e., for non-floating domains}.
     \hspace{.1in}(c)
\end{aligned}\right.
\end{equation}

{\em The system \eqref{ddsystemsh} is not a discretization of the system \eqref{ddsystems}} because  $\mathcal{A}_n^h(\cdot,\cdot)\ne \mathcal{A}_n(\cdot,\cdot)$, i.e., the former is defined with respect to $\Omega_n^h\cup\Gamma_n^h\cup\widehat\Gamma_n^h$ whereas the latter is defined in terms of $\Omega_n\cup\Gamma_n\cup\widehat\Gamma_n$. However, this observation is unimportant because what is true is that
\begin{equation}\label{goalgoaldis}
\begin{aligned}
&\mbox{\em the {global discrete solution $u^{dd,h}(\xb)$}, i.e. the solution such that}
\\ 
& u^h_n(\xb)=u^{dd,h}(\xb)|_{\Omega_n^h\cup\Gamma_n^h\cup\widehat\Gamma_n^h}, \mbox{\em is the same as the solution  $u^h(\xb)$}
\\
& \mbox{\em of the single-domain FE system \eqref{eq:weak_nonlochh}}
\end{aligned}
\end{equation}
which is, after all, the goal \eqref{goalgoal} we want to achieve. The truthfulness of \eqref{goalgoaldis} is verified following the same steps as those used to prove Proposition \ref{prop-equivalent} with, of course, $\Omega_n\cup\Gamma_n\cup\widehat\Gamma_n$ replaced by $\Omega_n^h\cup\Gamma_n^h\cup\widehat\Gamma_n^h$. 

When we define, in Section \ref{matrixform}, the matrix form of \eqref{ddsystemsh}, it is useful to  differentiate between the bilinear forms in the two cases in (\ref{ddsystemsh}a). First, because $\Gamma_n=\emptyset$ for {\em floating domains}, we have 
\begin{equation}\label{aflh}
\begin{aligned}
 & \mathcal{A}_n^{h}(u_n^h,v_n^h)=\mathcal{A}_n^{fl,h}(u_n^h,v_n^h) = 
 \\& \quad\int_{\Omega_n^h\cup\widehat\Gamma_n^h}
  \int_{\Omega_n^h\cup\widehat\Gamma_n^h}  
  \zeta_{\mcA}^h(\xb,\yb)^{-1} \big(v_n^h(\yb) - v_n^h(\xb)\big)\big(u_n^h(\yb) - u_n^h(\xb)\big)\gamma(\xb,\yb)d\yb  d\xb.
\end{aligned}
\end{equation}
For {\em non-floating domains}, the bilinear form involves integrals with respect to $\Gamma_n^h$. For such integrals, we have that either $v_n^h(\cdot)=0$ [because $v_n^h\in W_n^{0,h}$] or $u_n^h(\cdot)=g^h(\cdot)$ [because of the constraint (\ref{ddsystemsh}c)], so that then
\begin{equation}\label{anflh1}
  \mathcal{A}_n^{h}(u_n^h,v_n^h)={\mathcal A}^{n\!fl,h}_n(u_n^h,v_n^h) - {\mathcal A}^{g,h}_n(g^h,v_n^h),
\end{equation} 
where
\begin{equation}\label{anflh}
\begin{aligned}
&{\mathcal A}^{n\!fl,h}_n(u_n^h,v_n^h) =
\\ 
  &\qquad 
  \int_{\Omega_n^h\cup\widehat\Gamma_n^h}
  \hspace{-.02in}  \int_{\Omega_n^h\cup\widehat\Gamma_n^h}   
  \zeta_{\mcA}^h(\xb,\yb)^{-1} \big(v_n^h(\yb) - v_n^h(\xb)\big)\big(u_n^h(\yb) - u_n^h(\xb)\big)\gamma(\xb,\yb)d\yb  d\xb
  \\&\qquad\qquad+
{2\int_{\Omega_n^h\cup\widehat\Gamma_n^h}u_n^h(\xb)v_n^h(\xb)\int_{\Gamma_n^h}
\zeta^h_{\mcA}(\xb,\yb)^{-1}
\gamma(\xb,\yb) d\yb  d\xb}
\end{aligned}
\end{equation}
and
\begin{equation}\label{anflhg}
\begin{aligned}
{\mathcal A}^{g,h}_n(g^h,v_n^h)= 
  &\int_{\Omega_n^h\cup\widehat\Gamma_n^h}^{}\int_{\Gamma_n^h}{}\zeta_{\mcA}^h(\xb,\yb)^{-1}
g^h(\yb)v_n^h(\xb)
\gamma(\xb,\yb) d\yb  d\xb  
\\&\qquad
+\int_{\Gamma_n^h}{}\int_{\Omega_n^h\cup\widehat\Gamma_n^h}{}  \zeta^h_{\mcA}(\xb,\yb)^{-1}g^h(\xb)v_n^h(\yb) 
 \gamma(\xb,\yb)  d\yb  d\xb.
  \end{aligned}
\end{equation}
Note that in \eqref{aflh}, \eqref{anflh}, and \eqref{anflhg}, both $u_n^h(\cdot)$ and $v_n^h(\cdot)$ are evaluated only at points in $\Omega_n^h\cup\widehat\Gamma_n^h$ and $g_n^h(\cdot)$ is evaluated only at points in $\Gamma_n^h$.

In light of \eqref{aflh}, \eqref{anflh}, and \eqref{anflhg},
\eqref{ddsystemsh} can be rewritten as, for {\em floating domains},
\begin{equation}\label{fleq}
\left\{\begin{aligned}
\mathcal{A}_n^{fl,h}(u_n^h,v_n^h) &= \mcF_n^h(v_n^h)
\hspace{2.9in}(a)
\\
 u_n^h(\xb) &= u_{n'}^h(\xb) \,\,
    \quad \mbox{$\forall\,\xb\in \widehat\Gamma_{n}^h \cap \widehat{\Gamma}_{n'}^h$\,\,  \em for $n'=n+1,\ldots,N_s$}
    \qquad\quad\,\,\,(b)
\end{aligned}\right.
\end{equation}
 and for {\em non-floating domains}
\begin{equation}\label{nfleq}
\left\{\begin{aligned}
\mathcal{A}_n^{n\!fl,h}(u_n^h,v_n^h) &= \mcF_n^h(v_n^h) + {\mathcal A}^{g,h}_n(g^h,v_n^h)
\hspace{1.99in}(a)
\\
 u_n^h(\xb) &= u_{n'}^h(\xb) \,\,
     \quad \mbox{$\forall\,\xb\in \widehat\Gamma_{n}^h \cap \widehat{\Gamma}_{n'}^h$\,\,  \em for $n'=n+1,\ldots,N_s$.}
\qquad\quad(b)
\end{aligned}\right.
\end{equation}

\begin{remark}\label{rem5555}
{\em
For any domain indexed by $n$, each of the subdomain problems in \eqref{fleq} and \eqref{nfleq} is coupled, through (\ref{fleq})b or (\ref{nfleq})b, to other domains indexed by $n'$ with $n'\ne n$. Of course, this defeats the goal of domain decomposition which is to construct {\em uncoupled} subdomain problems so that, e.g., parallelization can be realized. This becomes the task for algorithms of obtaining {\em solutions} of the subdomain problems. Further comments in this regard are provided in Section \ref{sec:conclusion}.
}\qquad$\Box$
\end{remark}

\begin{remark}\label{addingc}
{\em 
If instead of (\ref{eq:strong_nonloc1}a) we consider
$$  
 -2\int_{\widehat\Omega\cup\Gamma\cup\Gamma_{N\!eumann}} 
\big(u(\yb) -u(\xb)\big)\gamma(\xb,\yb)d\yb 
+c(\xb)u(\xb)
= f_{\widehat\Omega}(\xb)
$$
with $c(\xb)>0$ so that the bilinear form ${\mathcal A}(u,v)$ in \eqref{eq:bil1} has the additional term $\int_\Omega c(\xb)u(\xb)v(\xb)d\xb$, then that bilinear form is coercive even for floating subdomains. In this case, the design of solution methods for \eqref{ddsystemsh}, and in particular for \eqref{fleq}, becomes substantially simpler.
}\qquad$\Box$
\end{remark}

\begin{remark}\label{remark888}
{\em The discussion that includes \eqref{aflh}--\eqref{nfleq} as well as the comments made in Remarks \ref{rem5555} and \ref{addingc} about the discrete bilinear form ${\mathcal A}_n^h(\cdot,\cdot)$ and the discrete subdomain system \eqref{ddsystemsh} also hold for the continuous bilinear form ${\mathcal A}_n(\cdot,\cdot)$ defined in \eqref{eq:bic} and continuous subdomain system \eqref{ddsystems}. 
}\qquad$\Box$
\end{remark}

%%%%%%%%%%%%%%
\subsubsection{Matrix form of the multi-domain finite element system}\label{matrixform}

Based on the numbering introduced in Section \ref{sec1-1}, we let

\hangtwo $X_{\Omega\cup\Gamma}=\{\xb_i\}_{i=1}^{\widetilde N^h}$ denote the set of nodes in the grid used for the FE discretization of the single-domain $\Omega\cup\Gamma$

\hangtwo $X_\Omega=\{\xb_i\}_{i=1}^{N^h}$ 
and
$X_\Gamma=\{\xb_i\}_{i=N^h+1}^{\widetilde N^h}$ denote the set of nodes in $\Omega$ and  $\Gamma$, respectively. 

\noindent The sets $X_\Omega\subset X_{\Omega\cup\Gamma}$ and $X_\Gamma\subset X_{\Omega\cup\Gamma}$ are disjoint and $X_\Omega\cup X_\Gamma=X_{\Omega\cup\Gamma}$. For $n=1,\ldots,N_s$, let 

\hangtwo $\widetilde X_n^h=\{\xb_i^n\}_{i=1}^{\widetilde N_n^h}$ denote a local numbering of the set of nodes in $\Omega_n^h\cup\widehat\Gamma_n^h\cup\Gamma_n^h$; clearly, by construction, $\widetilde X_n^h\subset X_{\Omega\cup\Gamma}$

\hangtwo $X_n^h=\{\xb_i^n\}_{i=1}^{N_n^h}$ and  $X_{\Gamma_n^h}^h=\{\xb_i^n\}_{N_n^h+1}^{\widetilde N_n^h}$ denote the nodes in $\widetilde X_n^h$ located in $\Omega_n^h\cup\widehat\Gamma_n^h$ and $\Gamma_n^h$, respectively.

\noindent The sets $X_n^h\subset \widetilde X_n^h$ and $X_{\Gamma_n^h}^h\subset \widetilde X_n^h$ are disjoint and $X_n^h\cup X_{\Gamma_n^h}^h=\widetilde X_n^h$. Note that for $n'\ne n$, the sets $X_n^h$ and $X_{n'}^h$ overlap whenever $\widehat\Gamma_n^h\cap\widehat\Gamma_{n'}^h\ne\emptyset$ and similarly for the sets $X_{\Gamma_n^h}$ and $X_{\Gamma_{n'}^h}$. Also note that if $\Omega_n^h$ is a floating domain, then the set $X_{\Gamma_n}^h$ is vacuous.

\medskip\noindent
\underline{\bf Bilinear forms in matrix notation}. We first consider the conversion of (\ref{fleq}a) and (\ref{nfleq}a) to matrix notation. Corresponding to the nodes in $X_{\Omega\cup\Gamma}$, we have the set of basis functions $\{\phi_i(\xb)\}_{i=1}^{\widetilde N^h}$ whose span is used to define the finite element space $W^h$ for the single-domain finite element system. We introduce the set of basis functions corresponding to each subdomain. For $n=1,\ldots,N_s$, let

\hangzero $\{\phi_i^n(\xb)\}_{i=1}^{\widetilde N_n^h}$ denote a local numbering of the basis functions in the spanning set for $W^h$ which correspond to the nodes in $\widetilde X_n^h$. 

\noindent By construction, $\{\phi_i^n(\xb)\}_{i=1}^{\widetilde N_n^h}$ spans the finite element space $W_n^h$. We then 

\hangzero divide $\{\phi_i^n(\xb)\}_{i=1}^{\widetilde N_n^h}$ into the sets $\{\phi_i^n(\xb)\}_{i=1}^{N_n^h}$ and $\{\phi_i^n(\xb)\}_{i=N_n^h+1}^{\widetilde N_n^h}$ that correspond to nodes in $X_n^h$ and $X_{\Gamma_n^h}^h$, respectively.

\noindent Note that if $\Omega_n$ is a floating domain, then $ N_n^h=\widetilde N_n^h$ so that the set $\{\phi_i^n(\xb)\}_{i=N_n^h+1}^{\widetilde N_n^h}$ is vacuous.

Let $\vecun$ denote an $N^h_n$-vector of nodal values of a function $u^h_n(\xb)$ defined for the nodes in $X_n^h$ (i.e., nodes in $\Omega^h_n\cup\widehat\Gamma^h_n$) and let $\vecgn$ denote the $(\widetilde N^h_n-N^h_n)$-vector of nodal values of $g^h(\xb)$ defined for the nodes in $X_{\Gamma_n^h}^h$ (i.e., nodes in $\Gamma^h_n$). Then, we have that
\begin{equation}\label{unhphi}
    u^h_n(\xb) = 
    \left\{\begin{tabular}{ll}
    $\displaystyle \sum_{j=1}^{N^h_n} \vecupn_j \phi_j^n(\xb)
    + \sum_{j=N^h_n+1}^{\widetilde N^h_n} \vecgpn_j \phi_j^n(\xb)$
    & {\em for non-floating domains}
    \\
    $\displaystyle\sum_{j=1}^{N^h_n} \vecupn_j \phi_j^n(\xb)$
     &{\em for floating domains.}
    \end{tabular}
    \right.
\end{equation}
Note that components of $\vecun$ are ordered according the local indexing of nodes.
We then define, for $n=1,\ldots,N_s$, the entries of the $N_n^h\times N_n^h$ matrix $\mathbb{A}_n$ as
\begin{equation}\label{anij}
  (\mathbb{A}_n)_{ij} =
   \left\{\begin{tabular}{ll}
   \hspace{-.1in}${\mathcal A}^{n\!fl,h}_n(\phi_j^n,\phi_i^n)$
   & {\em for non-floating domains}
   \\[1.5ex]
   \hspace{-.1in}${\mathcal A}^{fl,h}_n(\phi_j^n,\phi_i^n)$
   & {\em for floating domains}
    \end{tabular}
    \right.
    \mbox{for $i,j=1,\ldots,N_n^h$}
\end{equation}
and, for $i=1,\ldots,N_n^h$, the components of the $N_n^h$-vector $\vecbn$ as
\begin{equation}\label{bni}
   \vecbpn_i =
   \left\{\begin{tabular}{ll}
   \hspace{-.1in}${\mathcal F}^{h}_n(\phi_i^n)$
  $\displaystyle -\sum_{j=N^h+1}^{\widetilde N^h}{\mathcal A}^{g,h}_n(\phi_j^n,\phi_i^n)\vecgpn_j$
   & {\em for non-floating domains}
   \\[1.5ex]
   \hspace{-.1in}${\mathcal F}^{h}_n(\phi_i^n)$
   & {\em for floating domains}.
    \end{tabular}
    \right.
\end{equation}
Then, the finite element problems (\ref{fleq}a) and (\ref{nfleq}a) both have the matrix notation equivalent
\begin{equation}\label{aueqb}
       {\mathbb A}_n\vecun = \vecbn\qquad\mbox{for $n=1,\ldots,N_s$}.
\end{equation}

\medskip\noindent
\underline{\bf Constraints in matrix notation}. We next turn to the conversion of the constraints in (\ref{ddsystemsh}b) [or equivalently (\ref{fleq}b) and (\ref{nfleq}b)] to matrix notation. 

We first transform those constraints to vector notation. We have global and local indices of nodes. Thus, if the node $\xb_j^n\in X_n$ corresponds to the node $\xb_i\in X$, i.e., if we have $\xb_i = \xb_j^n$, then $i$ and $j$ are the global and local indices, respectively, for the same node. We define a mapping from global to local indices, specifically, for a globally indexed node $\xb_i\in\widehat\Gamma_n^h$, we let
$$
I_{ni}= \mbox{\em local index of the node $\xb_i\in \widehat\Gamma_n^h$.}
$$
Now suppose that $\widehat\Gamma_n^h\cap\widehat\Gamma_{n'}^h\ne\emptyset$. Then, for a globally indexed node $\xb_i\in \widehat\Gamma_n^h\cap\widehat\Gamma_{n'}^h$, we have that
$$
   \xb_i = \xb_j^n = \xb_{j'}^{n'}
   \qquad\mbox{where $j=I_{ni}$ and $j'=I_{n'i}$}.
$$
Because the FE approximation $u_n^h(\xb)$ is uniquely determined by its nodal values, the constraints in (\ref{ddsystemsh}b) can be equivalently expressed as
\begin{equation}\label{ieqip}
     \vecupn_j = (\vecu_{n'})_{j'}\quad
     \mbox{for all nodes \quad$\xb_i\in \widehat\Gamma_{n}^h \cap \widehat{\Gamma}_{n'}^h$},\quad
    n'= n+1,\ldots,N_s.
\end{equation}
We keep in mind that, as was the case for the continuous multi-domain system (see Remark \ref{redunconst}), the constraints in \eqref{ieqip} are not independent. For example, if $\widehat\Gamma_{1}^h \cap \widehat{\Gamma}_{2}^h\cap \widehat{\Gamma}_{3}^h\ne\emptyset$, then, for a node $\xb_i$ in that domain, we have from \eqref{ieqip} that
$(\vecu_{1})_{j} = (\vecu_{2})_{j'}$, $(\vecu_{1})_{j} = (\vecu_{3})_{j''}$, and $(\vecu_{2})_{j'} = (\vecu_{3})_{j''}$, where $j=I_{1i}$, $j'=I_{2i}$, and $j''=I_{3i}$.
Clearly, these three equations are not independent.

Our task is then reduced to expressing the constraints in \eqref{ieqip} in an economical matrix form. Here, we mimic the process given in \cite{mathew2008domain} for DD in the local PDE case, In fact, we construct an $M\times \sum_{n=1}^{N_s}{N_n^h}$ matrix ${\mathbb M}$, where $M$ is specified below, of the form
$$
   {\mathbb M} = ( {\mathbb M}_1 \cdots {\mathbb M}_{N_s} )
$$
such that the constraints in \eqref{ieqip} can be equivalently expressed as
\begin{equation}\label{matrixcoup}
\sum_{n=1}^{N_s}  {\mathbb M}_n \vecun = 0,
\end{equation}
where ${\mathbb M}_n$, $n=1,\ldots,N_s$, are $M\times N_n^h$ matrices. We construct two such matrices, one for the constraints in \eqref{ieqip}, the other for an equivalent non-redundant set of constraints. To this end, letting $\widehat\Gamma^h=\cup_{n=1}^{N_s}\widehat\Gamma_n^h$,
\begin{equation}\label{degree}
\begin{aligned}
&\mbox{\em for each node in $\xb_i\in\widehat\Gamma^h$, we define the set}
\\&
\qquad\qquad \theta(\xb_i)= \{n\,\,:\,\, \xb_i\in\widehat\Gamma_n^h\}
\\&
\mbox{\em and let\,\,\, $m(\xb_i)=$ cardinality of the set $\theta(\xb_i)$}
\end{aligned}
\end{equation}
so that $\theta(\xb_i)$ consists of the indices of all the subdomains $\widehat\Gamma_n^h$ that contain the globally indexed node $\xb_i$ and $m(\xb_i)$ denotes the number of distinct subdomains which the node $\xb_i$ belongs to.

\medskip\noindent
\underline{\em Non-redundant constraints -- full-rank matrix ${\mathbb M}$}.
For each node $\xb_i\in\widehat\Gamma^h$,

\qquad - arrange the indices in $\theta(\xb_i)$ in increasing order

\qquad - for each {\em consecutive} pair of indices, impose one constraint.

\noindent This results in, with $N_{\widehat\Gamma^h}$ denoting the number of nodes in $\widehat\Gamma^h$,
$$
     M = \sum_{i=1}^{N_{\widehat\Gamma^h}} \big( m(\xb_i) - 1 \big)
$$
non-redundant constraints. For example, if $\widehat\Gamma_{1}^h \cap \widehat{\Gamma}_{2}^h\cap \widehat{\Gamma}_{3}^h\ne\emptyset$, then, for a node $\xb_i$ in that domain, we now have that
$(\vecu_{1})_{j} = (\vecu_{2})_{j'}$ and $(\vecu_{2})_{j'} = (\vecu_{3})_{j''}$, where $j=I_{1i}$, $j'=I_{2i}$, and $j''=I_{3i}$. Clearly, these two equations are not redundant and together imply the constraint $(\vecu_{1})_{j}  = (\vecu_{3})_{j''}$ from \eqref{ieqip} that is now missing.

The entries of the matrix ${\mathbb M}$ can be determined as follows: set $k=0$ and then,

\hangtwo\qquad - for $i=1,\ldots, N_{\widehat\Gamma^h}$

\hangtwo\qquad - for each pair {$(n,n')$}, $n<n'$, of consecutive indices in $\theta(\xb_i)$

\hangtwo\qquad - set $k\longleftarrow k+1$

\hangtwo\qquad - for $j$ and $j'$ such that $\xb^n_{j}=\xb^{n'}_{j'}=\xb_i$, 
$$
\mbox{{set}} \,\,\,\, ({\mathbb M}_n)_{kj} = 1  \qquad
({\mathbb M}_{n'})_{kj'} = -1 
$$
and all other entries in the $k$-th row of ${\mathbb M}$ to zero.

\medskip\noindent
\underline{\em Redundant constraints}
We proceed as we do above for the non-redundant set of constraints, except that now we do not require that $n<n'$ be consecutive indices in $\theta(\xb_i)$, i.e., we impose a constraint for every distinct pair $n<n'$ of indices in $\theta(\xb_i)$. This approach lends itself better for parallelization compared to the use of non-redundant constraints; see, e.g., \cite{mathew2008domain} for a discussion in the local DD setting. Thus, we impose all the constraints in \eqref{ieqip}. Note that redundant constraints are caused only for $m(\xb_i)\ge 3$. {Because there are $m(\xb_i)$ distinct indices in $\theta(\xb_i)$}, we have that the number of rows in the matrix ${\mathbb M}$ is now given by
$$
M = \sum_{i=1}^{N_{\widehat\Gamma^h}} \frac12 m(\xb_i)\big( m(\xb_i) - 1 \big).
$$

%%%%%%%%%%%%%%%%%%%%%%%%%%%%%%%%%%%%%%%%%%%%%%%%%%%%
%%%%%%%%%%%%%%%%%%%%%%%%%%%%%%%%%%%%%%%%%%%%%%%%%%%%
\section{Concluding remarks}\label{sec:conclusion}

We have defined and analyzed a general framework for the construction of domain decomposition methods for nonlocal problems that achieves the goal stated in \eqref{goalgoal2} or, more precisely, in \eqref{goalgoal}. However, there is still work to be done because we have not met a second goal which is that the nonlocal DD method is amenable to parallelization. What we have so far are $N_s$ systems \eqref{fleq}--\eqref{nfleq} in which the subdomain systems are coupled through the constraints (\ref{fleq}b) and (\ref{nfleq}b) or, equivalently, the subdomain matrix systems in \eqref{aueqb} that are coupled through the constraints in \eqref{matrixcoup}. Such couplings prevent the direct use of \eqref{fleq}--\eqref{nfleq} (or equivalently \eqref{aueqb}--\eqref{matrixcoup}) for achieving the second goal.

At a similar stage in the development DD algorithms in the local PDE setting, one is faced with the analogous situation; for example, in the non-overlapping DD setting, there is coupling of the subdomain problems, say in the matrix formulation, at the nodes located along the common boundaries between subdomains. The uncoupling between subdomain problems is then effected through the design of {\em solution methods} in which the bulk of the computational effort is borne by steps in those methods that are parallelizable and for which the non-paralellizable steps and the communications between subdomain problems, i.e., between processors, is kept to a minimum. 

In the nonlocal DD setting, solution methods have to be designed to meet the same criterion: the bulk of the computational costs has to be borne by parallelizable steps. In a follow-up paper, we will develop, analyze, and implement such methods. Parallelizable solution methods for local PDE non-overlapping DD will be generalized to the nonlocal setting. For example, Lagrange multiplier methods, e.g., FETI \cite{farhat1991method}, Arlequin methods \cite{BenDhia_05_IJNME}, and optimization-based DD methods \cite{Gunzburger_00_AMC,Gunzburger_99_CMA,Bochev_17a_SISC}, all of which are in use for local PDE non-overlapping DD and which are also all good candidates for generalization to the nonlocal DD setting. Other local DD solution methods could also be considered for generalization. As is the case for the framework developed in this paper, generalizations of solution methods will pose challenges because of nonlocality. One thing to keep in mind is that a preferred solution method in the local setting may or may not remain so when generalized to the nonlocal setting.

%%%%%%%%%%%%%%%%%%%%%%%%%%%%%%%%%%%%%%%%%%%%%%%%%%%%
%%%%%%%%%%%%%%%%%%%%%%%%%%%%%%%%%%%%%%%%%%%%%%%%%%%%
\section*{Acknowledgments}
Sandia National Laboratories is a multimission laboratory managed and operated by National Technology and Engineering Solutions of Sandia, LLC., a wholly owned subsidiary of Honeywell International, Inc., for the U.S. Department of Energy's National Nuclear Security Administration under contract DE-NA-0003525. This paper describes objective technical results and analysis. Any subjective views or opinions that might be expressed in this paper (SAND2020-8735) do not necessarily represent the views of the U.S. Department of Energy or the United States Government.

This work was supported by the Sandia National Laboratories (SNL) Laboratory-directed Research and Development (LDRD) program, and the U.S. Department of Energy, Office of Science, Office of Advanced Scientific Computing Research under Award Number DE-SC-0000230927 and under the Collaboratory on Mathematics and Physics-Informed Learning Machines for Multiscale and Multiphysics Problems (PhILMs) project.

\bibliography{domdecomp.bib}

%%%%%%%%%%%%%%%%%%%%%%%%%%%%%%%%%%%%%%%%%%%%%%%%%%%%
%%%%%%%%%%%%%%%%%%%%%%%%%%%%%%%%%%%%%%%%%%%%%%%%%%%%
\appendix
%%%%%%%%%%%%%%%%%%%%%%%%%%%%%%%%%%%%%%%%%%%%%%%%%%%%
%%%%%%%%%%%%%%%%%%%%%%%%%%%%%%%%%%%%%%%%%%%%%%%%%%%%
\section{Proof of Lemma \ref{lemma111}}\label{main-proof}
It is convenient for what follows to introduce, for each $n=1,\ldots,N_s$ such that $\Gamma_n\ne\emptyset$, the splitting
\begin{equation*}\label{gammasplit}
\Gamma_n = \Gamma_n^*\cup\Gamma_n^{\dag}
\quad\mbox{with}\quad 
\left\{\begin{aligned}
& 
  \Gamma_n^*\subset\Gamma_n
  \quad\mbox{such that}\quad
  \Gamma_n^*\cap(\cup_{n'=1,\,n'\ne n}^{N_s}\Gamma_{n'})=\emptyset
\\
&\Gamma_n^{\dag}=\Gamma_n\setminus\Gamma_n^*\subset\Gamma_n
\end{aligned}\right.
\end{equation*}
so that $\Gamma_n^*$ (resp. $\Gamma_n^{\dag}$) are the disjoint parts of $\Gamma_n$ that do not (resp. do) overlap with any other $\Gamma_{n'}$ with $n'\ne n$. The blue regions in Figure \ref{fig2.2}-right illustrate examples of the sets $\Gamma_n^\dag$. Note that, by definition,
$\Gamma_n^*\cap\Gamma_n^{\dag}=\emptyset$. 
With these definitions in hand, we have that, for $n=1,\ldots,N_s$, 
$
\Omega_n\cup\widehat\Gamma_n\cup\Gamma_n =\Omega_n\cup\widehat\Gamma_n\cup\Gamma_n^*\cup\Gamma_n^{\dag}
$.
We can then express \eqref{eq:zeta} as
\begin{equation*}\label{eq:zeta2}
\zeta_{\mathcal{A}}(\xb,\yb)
=
\sum_{n=1}^{N_s}
\mcX_{\Omega_n\cup\Gamma_n^*\cup\widehat\Gamma_n\cup\Gamma_n^{\dag}}(\xb)
\mcX_{\Omega_n\cup\Gamma_n^*\cup\widehat\Gamma_n\cup\Gamma_n^{\dag}}(\yb).
\end{equation*}
Note that among the sets $\Omega_n,\;\Gamma_n^*,\; \widehat\Gamma_n$, and $\Gamma_n^\dag$ the only two that may possibly overlap with other sets are $\Gamma_n^{\dag}$ and $\widehat{\Gamma}_n$. For this reason and for ease of notation, we further introduce the set $\widetilde\Gamma_n=\Gamma_n^{\dag}\cup\widehat{\Gamma}_n$. Based on this consideration, we split the outer integral of the bilinear form $\mcA_n$ into a set that does not overlap with any other sets (i.e. $\Omega_n\cup\Gamma_n^*$) and $\widetilde\Gamma_n$. We have:
\begin{equation*}
 \begin{aligned}\label{eq:bil}
  &\mathcal{A}_n(u_n,v_n) \\
  &\qquad = \int_{\Omega_n\cup\Gamma_n^*\cup\widetilde\Gamma_n}  \int_{\Omega_n\cup\Gamma_n^*\cup\widetilde\Gamma_n} \zeta_{\mathcal{A}}(\mathbf{x},\mathbf{y})^{-1}\big(v_n(\mathbf{y}) - v_n(\mathbf{x})\big)\big(u_n(\mathbf{y}) - u_n(\mathbf{x})\big)\gamma(\mathbf{x},\mathbf{y})d\mathbf{y}  d\mathbf{x}\\[2mm]
  &\qquad = \mathcal{A}_n^{disjoint}(u_n,v_n) + \mathcal{A}_n^{overlap}(u_n,v_n),
 \end{aligned}
\end{equation*}
where
\begin{equation*}
 \begin{aligned}
  &\mathcal{A}_n^{disjoint}(u_n,v_n) := \\
  &\qquad \int_{\Omega_n\cup\Gamma_n^*}  \int_{\Omega_n\cup\Gamma_n^*\cup\widetilde\Gamma_n} \zeta_{\mathcal{A}}(\mathbf{x},\mathbf{y})^{-1}\big(v_n(\mathbf{y}) - v_n(\mathbf{x})\big)\big(u_n(\mathbf{y}) - u_n(\mathbf{x})\big)\gamma(\mathbf{x},\mathbf{y})d\mathbf{y}  d\mathbf{x},
 \end{aligned}
\end{equation*}
\begin{equation*}
 \begin{aligned}
  &\mathcal{A}_n^{overlap}(u_n,v_n) := \\
  &\qquad \int_{\widetilde\Gamma_n}  \int_{\Omega_n\cup\Gamma_n^*\cup\widetilde\Gamma_n} \zeta_{\mathcal{A}}(\mathbf{x},\mathbf{y})^{-1}\big(v_n(\mathbf{y}) - v_n(\mathbf{x})\big)\big(u_n(\mathbf{y}) - u_n(\mathbf{x})\big)\gamma(\mathbf{x},\mathbf{y})d\mathbf{y}  d\mathbf{x}.
 \end{aligned}
\end{equation*}
To simplify the notation, we let $w_n(\xb,\yb)= (v_n(\mathbf{y}) - v_n(\mathbf{x})\big)(u_n(\mathbf{y}) - u_n(\mathbf{x})\big)\gamma(\mathbf{x},\mathbf{y})$ and $w(\xb,\yb) = (v(\mathbf{y}) - v(\mathbf{x})\big)(u(\mathbf{y}) - u(\mathbf{x})\big)\gamma(\mathbf{x},\mathbf{y})$. Note that whenever $\xb$ or $\yb$ belong to a set that overlaps with other sets, {the definitions of $u_n$ and $v_n$ in \eqref{uequn} guarantee that} $w_n(\xb,\yb)=w(\xb,\yb)$. We first analyze $\mcA_n^{disjoint}$; we have that
\begin{equation*}
\begin{aligned}
\mathcal{A}_n^{disjoint}(u_n,v_n) 
&=\int_{\Omega_n\cup\Gamma_n^*}  \int_{\Omega_n\cup\Gamma_n^*\cup\widetilde\Gamma_n} w(\mathbf{x},\mathbf{y})d\yb d\xb\\
&= \int_{\Omega_n\cup\Gamma_n^*}  \int_{\Omega_n\cup\Gamma_n^*\cup\widetilde\Gamma_n} w(\mathbf{x},\mathbf{y})d\yb d\xb\\
&+  \int_{\Omega_n\cup\Gamma_n^*} \int_{(\Omega\cup\Gamma)\setminus(\Omega_n\cup\Gamma_n^*\cup\widetilde\Gamma_n)}w(\mathbf{x},\mathbf{y})d\mathbf{y}  d\mathbf{x}  
 \end{aligned}
\end{equation*}
where the first equality follows from the fact that $\zeta_{\mathcal{A}}(\mathbf{x},\mathbf{y})=1$ for $\mathbf{x} \in  \Omega_n \cup \Gamma_n^*$ and $\yb\in\Omega_n\cup\Gamma_n^*\cup\widetilde\Gamma_n$ and {the second inequality from the fact that $\gamma(\mathbf{x},\mathbf{y})=0$ for $\mathbf{x} \in  \Omega_n \cup \Gamma_n^*$ and $\yb\in(\Omega\cup\Gamma)\setminus(\Omega_n\cup\Gamma_n^*\cup\widetilde\Gamma_n)$.} Hence,
\begin{equation}\label{disjoint}
 \begin{aligned}
  \sum\limits_{n=1}^{N_s}\mathcal{A}_n^{disjoint}(u_n,v_n) = \int_{\bigcup\limits_{n=1}^{N_s}(\Omega_n\cup\Gamma_n^*)}  \int_{\Omega\cup\Gamma} w(\mathbf{x},\mathbf{y})d\mathbf{y}  d\mathbf{x}.
 \end{aligned}
\end{equation}
By definition of $\zeta_{\mathcal{A}}$, for any $\mathbf{y} \in (\Omega_n\cup\Gamma_n^*)$ and for any $\mathbf{x} \in \widetilde\Gamma_n$ we have $\zeta_{\mathcal{A}}(\mathbf{x},\mathbf{y})=1$, hence 
\begin{equation*}
 \begin{aligned}
  \mathcal{A}_n^{overlap}(u_n,v_n) &= \int_{\widetilde\Gamma_n}  \int_{\Omega_n\cup\Gamma_n^*} w_n(\mathbf{x},\mathbf{y})\, d\mathbf{y}  d\mathbf{x}\\
  &+ \int_{\widetilde\Gamma_n}  \int_{\widetilde\Gamma_n} \zeta_{\mathcal{A}}(\mathbf{x},\mathbf{y})^{-1}
  w_n(\mathbf{x},\mathbf{y})\,d\mathbf{y}  d\mathbf{x}\\
  &= \mathcal{A}_n^{overlap,I}(u_n,v_n) + \mathcal{A}_n^{overlap,II}(u_n,v_n).  
 \end{aligned}
\end{equation*}
We introduce the {partition} $\{\Lambda_i\}_{i=1}^{N_\lambda}$ of $\bigcup_{n=1}^{N_s}\widetilde\Gamma_n$, i.e. 
$$\bigcup_{i=1}^{N_\lambda}\Lambda_i = \bigcup_{i=1}^{N_s}\widetilde\Gamma_n, \quad {\rm and} \quad \Lambda_i\cap\Lambda_j=\emptyset, \;{\rm for} \; i\neq j,
$$
such that there exists an index set $I_n$, $n=1, \ldots N_s$, for which 
$$
\bigcup_{i\in I_n}\Lambda_i = \widetilde\Gamma_n.
$$
{In practice, this means that for every $n=1,\ldots N_s$ there exists a subset of the partition that provides a partition of $\widetilde\Gamma_n$. Note that such a partition always exists.} Then we have
\begin{equation*}
\begin{aligned}
      \sum\limits_{n=1}^{N_s}\mathcal{A}_n^{overlap,I}(u_n,v_n) 
&= \sum\limits_{n=1}^{N_s} \int_{\widetilde\Gamma_n}  \int_{\Omega_n\cup\Gamma_n^{*}}w_n(\mathbf{x},\mathbf{y}) d\mathbf{y}  d\mathbf{x} \\
&= \sum\limits_{n=1}^{N_s}\sum\limits_{i\in I_n} \int_{\Lambda_i}  \int_{\Omega_n\cup\Gamma_n^{*}}w_n(\mathbf{x},\mathbf{y}) d\mathbf{y}  d\mathbf{x}\\
&= \sum\limits_{n=1}^{N_s}\sum\limits_{i=1}^{N_\lambda} \int_{\Lambda_i}\int_{\Omega_n\cup\Gamma_n^{*}}w_n(\mathbf{x},\mathbf{y}) d\mathbf{y}d\xb \\
&= \sum\limits_{i=1}^{N_\lambda}\int_{\Lambda_i}\left[\sum\limits_{n=1}^{N_s}  \int_{\Omega_n\cup\Gamma_n^{*}}w_n(\mathbf{x},\mathbf{y}) d\yb\right]d\xb \\
&=\int_{\bigcup_{n=1}^{N_s} \widetilde\Gamma_n}\int_{\bigcup_{n=1}^{N_s}(\Omega_n\cup\Gamma_n^{*})}w(\mathbf{x},\mathbf{y}) d\mathbf{y}  d\mathbf{x},
\end{aligned}
\end{equation*}
where the second equality follows from the fact that $\{\Lambda_i\}_{i\in I_n}$ form a disjoint covering of $\widetilde\Gamma_n$, the third from the fact that $\gamma(\xb,\yb)=0$ for all the extra terms, the fourth from the fact that the sums are independent and hence can be switched, and the fifth from the fact that the sets are disjoint.

\smallskip For $\mathcal{A}_n^{overlap,II}$ we have 
\begin{equation*}
\begin{aligned}
\sum\limits_{n=1}^{N_s}\mathcal{A}_n^{overlap,II}(u_n,v_n) 
&= \sum\limits_{n=1}^{N_s} \int_{\widetilde\Gamma_n}\int_{\widetilde\Gamma_n}\zeta_A(\xb,\yb)^{-1}w(\xb,\yb) d\mathbf{y}  d\mathbf{x}\\
&= \sum\limits_{n=1}^{N_s}\sum\limits_{i\in I_n}\sum\limits_{j\in I_n} \int_{\Lambda_i} \int_{\Lambda_j} \zeta_A(\xb,\yb)^{-1} w(\xb,\yb) d\mathbf{y}  d\mathbf{x},
\end{aligned}
\end{equation*}
where, again, the second equality follows from the fact that the sets $\{\Lambda_i\}_{i\in I_n}$ are disjoint and form a covering of $\widetilde\Gamma_n$.

Next, we introduce the index set $\mcI_{ij}(\xb,\yb)$ that contains all indices $n$ such that $(\xb,\yb)\in\Lambda_i\times\Lambda_j$ and $i,j\in I_n$. Formally,
\begin{equation}\label{eq:mcIn}
\mcI_{ij}(\xb,\yb)=\{n\in\{1,\ldots,N_s\} \; : \; \xb\in\Lambda_i, \; \yb\in\Lambda_j, \; {\rm and} \; i,\,j\in I_n\}.
\end{equation}
%
%\CV{[KV: The following could be deleted, since the notation $\mu_{ij}(\xb,\yb)$ is not used below] We let cardinality$(\mcI_{ij}(\xb,\yb))=\mu_{ij}(\xb,\yb)$; this is the number of times the pair $(\xb,\yb)\in\Lambda_i\times\Lambda_j$ is considered when solving all subproblems, i.e.
%
%$$
%{\mu_{ij}(\xb,\yb) = \sum_{n=1}^{N_s} \mcX(n \in \mcI_{ij}(\xb,\yb)).}
%$$
%}
Note that for $\xb\in \Lambda_i$ and $\yb\in\Lambda_j$ 
{\begin{equation}\label{eq:zetaNij}
\begin{aligned}
   \zeta_\mcA(\xb,\yb) 
&= \sum\limits_{n=1}^{N_s}\mcX_{\Omega_n\cup\Gamma_n^*\cup\widetilde\Gamma_n}(\xb)
   \mcX_{\Omega_n\cup\Gamma_n^*\cup\widetilde\Gamma_n}(\yb) \\
&= \sum\limits_{n=1}^{N_s} \mcX_{\mcI_{ij}(\xb,\yb)}(n)
   \mcX_{\Lambda_i}(\xb)\mcX_{\Lambda_j}(\yb)
=  \sum\limits_{n=1}^{N_s} \mcX_{\mcI_{ij}(\xb,\yb)}(n).
%\CV{= \mu_{ij}(\xb,\yb)}.
\end{aligned}
\end{equation}}
Thus,
\begin{equation*}
\begin{aligned}
\sum\limits_{n=1}^{N_s}&\mathcal{A}_n^{overlap,II}(u_n,v_n) 
= \sum\limits_{n=1}^{N_s}\sum\limits_{i\in I_n}\sum\limits_{j\in I_n} \int_{\Lambda_i} \int_{\Lambda_j} \zeta_A(\xb,\yb)^{-1} w(\xb,\yb) d\yb d\xb\\
&= \sum\limits_{n=1}^{N_s}\sum\limits_{i\in I_n}\sum\limits_{j\in I_n}\int_{\Lambda_i} \int_{\Lambda_j}
      \mcX_{\mcI_{ij}(\xb,\yb)}(n) \mcX_{\Lambda_i}(\xb)\mcX_{\Lambda_j}(\yb)\zeta_A(\xb,\yb)^{-1} w(\xb,\yb) d\yb d\xb\\
&= \sum\limits_{n=1}^{N_s}\sum\limits_{i=1}^{N_\lambda}\sum\limits_{j=1}^{N_\lambda} \int_{\Lambda_i} \int_{\Lambda_j} 
      \mcX_{\mcI_{ij}(\xb,\yb)}(n) \mcX_{\Lambda_i}(\xb)\mcX_{\Lambda_j}(\yb)\zeta_A(\xb,\yb)^{-1} w(\xb,\yb) d\yb d\xb\\
&= \sum\limits_{i=1}^{N_\lambda}\sum\limits_{j=1}^{N_\lambda}  \int_{\Lambda_i} \int_{\Lambda_j}\zeta_A(\xb,\yb)\zeta_A(\xb,\yb)^{-1}w(\xb,\yb) d\yb d\xb\\
& = \int_{\bigcup_{n=1}^{N^s}\widetilde\Gamma_n}\int_{\bigcup_{n=1}^{N^s}\widetilde\Gamma_n} w(\xb,\yb)d\yb d\xb.
\end{aligned}
\end{equation*}
Here, in the second equality we only added terms that are equal to 1 because all the indicator functions are active. This allows us to extend the sums over $i$ and $j$ in the third equality because all the extra terms are zero. Then, the fourth equality follows from \eqref{eq:zetaNij} and the fifth from the fact that all sets $\Lambda_i$ and $\Lambda_j$ are disjoint.

It then follows that 
\begin{equation}
\begin{aligned}
 &\sum\limits_{n=1}^{N_s}\mathcal{A}_n^{overlap}(u_n,v_n)  =  \sum\limits_{n=1}^{N_s}\Big(\mathcal{A}_n^{overlap,I}(u_n,v_n)  
+\mathcal{A}_n^{overlap,II}(u_n,v_n)\Big)\\
&=\int_{\bigcup_{n=1}^{N_s} \widetilde\Gamma_n}\Big( \int_{\bigcup_{n=1}^{N_s}(\Omega_n\cup\Gamma_n^{*})}w(\mathbf{x},\mathbf{y}) d\mathbf{y}  +\int_{\bigcup_{n=1}^{N_s} \widetilde\Gamma_n}w(\mathbf{x},\mathbf{y}) d\mathbf{y}  \Big)d\mathbf{x}\\
&= \int_{\bigcup_{n=1}^{N_s} \widetilde\Gamma_n} \int_{\Omega \cup \Gamma} w(\mathbf{x},\mathbf{y}) d\mathbf{y} d\mathbf{x}.
  \end{aligned}
  \end{equation}
  The last equality follows from the fact that $\bigcup_{n=1}^{N_s}(\Omega_n\cup\Gamma_n^{*})$ and $\bigcup_{n=1}^{N_s} (\Gamma_n^{\dag}\cup\widehat{\Gamma}_n)$ are two disjoint sets. For the same reason, we have
\begin{equation}
\begin{aligned}
& \sum\limits_{n=1}^{N_s} \mathcal{A}_n(u_n,v_n) =  \sum\limits_{n=1}^{N_s}\Big(\mathcal{A}_n^{disjoint}(u_n,v_n)  
+\mathcal{A}_n^{overlap}(u_n,v_n)\Big)\\
&\int_{\bigcup\limits_{n=1}^{N_s}(\Omega_n\cup\Gamma_n^*)}  \int_{\Omega\cup\Gamma} w(\mathbf{x},\mathbf{y})d\mathbf{y}  d\mathbf{x} + \int_{\bigcup_{n=1}^{N_s} \widetilde\Gamma_n} \int_{\Omega \cup \Gamma} w(\mathbf{x},\mathbf{y}) d\mathbf{y} d\mathbf{x}\\
&= \int_{\Omega \cup \Gamma} \int_{\Omega \cup \Gamma} w(\mathbf{x},\mathbf{y}) d\mathbf{y} d\mathbf{x} = \mathcal{A}(u,v).
\end{aligned}
\end{equation}

The proof for $\mcF$ directly follows from the same arguments and it is not reported.

\end{document}